\documentclass[10pt]{article}

\usepackage{mathptmx}       
\usepackage{helvet}         
\usepackage{courier}        
\usepackage{type1cm}        
\usepackage{makeidx}         
\usepackage{graphicx}        
\usepackage{multicol}        
\usepackage[bottom]{footmisc}

\makeindex             

\usepackage{amsmath}
\usepackage{amssymb}

\newtheorem{Th}{Theorem}
\newtheorem{Cor}{Corollary}
\newtheorem{Prop}{Proposition}
\newtheorem{Lm}{Lemma}
\newtheorem{Def}{Definition}

\newtheorem{Ques}{Question}
\newtheorem{rem}{Remark}

\newcommand{\ws}{{\widetilde{\*\S}}}
\newcommand{\Limpl}{\Longrightarrow}
\newcommand{\A}{{\cal A}}

\newcommand{\Z}{\mathbb{Z}}
\newcommand{\R}{\mathbb{R}}
\newcommand{\N}{\mathbb{N}}

\renewcommand{\S}{\mathbb{S}}

\newcommand{\cs}{{\cal S}}
\newcommand{\bs}{{\bf S}}

\renewcommand{\P}{{\cal P}}

\newcommand{\M}{{\cal M}}

\newcommand{\F}{{\cal F}}
\newcommand{\iN}{\*\N_{\infty}}

\newcommand{\e}{\varepsilon}
\renewcommand{\l}{\lambda}
\newcommand{\f}{\varphi}
\newcommand{\Mu}{{\rm M}}
\renewcommand{\a}{\alpha}
\renewcommand{\b}{\beta}
\renewcommand{\d}{\delta}
\newcommand{\D}{\Delta}

\newcommand{\s}{\sigma}

\renewcommand{\ss}{\subseteq}
\renewcommand{\l}{\lambda}

\newcommand{\la}{\langle}
\newcommand{\ra}{\rangle}

\newcommand{\range}{\mbox{range}}

\newcommand{\dom}{\mbox{dom}}
\renewcommand{\*}{\,^*\!}
\renewcommand{\o}{\,^\circ\,\!\!}

\newcommand{\all}{\forall}
\newcommand{\ex}{\exists}
\newcommand{\Liff}{\Longleftrightarrow}
\newcommand{\liff}{\longleftrightarrow}

\newcommand{\alls}{\all^{\mbox{st}}}
\newcommand{\exs}{\ex^{\mbox{st}}}
\newcommand{\st}{^{\mbox{st}}}
\newcommand{\il}{^{\mbox{int}}}
\usepackage{latexsym}
\title{Nonstandard analysis of the behavior of ergodic means of dynamical systems on very big finite probability spaces.}
\author{E.I. Gordon, L.Yu. Glebsky, C.W. Henson}

\begin{document}
\maketitle


\section{Introduction}

We discuss here the behavior of ergodic means of discrete time dynamical systems on a very big finite probability space $Y$ (discrete dynamical systems below). The G. Birkhoff Ergodic Theorem states the eventual stabilization of ergodic means of integrable functions for almost all points of the probability space. The trivial proof of this theorem for the case of finite probability spaces shows that this stabilization happens for those time intervals, whose length $n$ exceeds significantly the cardinality $|Y|$ of $Y$, i.e. $\frac n{|Y|}$ is a very big number.

For the case of a very big number $|Y|$ we introduce a huge class of functions on $Y$ including, for example all bounded functions, i.e. those functions, whose values are significantly less, than $|Y|$. 
Functions of this class are said to be $S$-integrable 
(the formula (\ref{S-int-1}) below). 
The class of $S$-integrable functions is an analog of the class of integrable 
functions on an infinite probability space.

We show that the behavior of ergodic means of $S$-integrable functions demonstrates some 
regularity even for those intervals, whose length is comparable with $|Y|$. The ergodic 
means $A_n$ and $A_m$ on the intervals of time $\{0,...,m-1\}=\bar m$ and 
$\{0,...,n-1\}=\bar n$ are approximately the same if 
$\frac n{|Y|}\approx\frac m{|Y|}\gg 0.$ It means that if we plot the points 
$(\frac n{|Y|}, A_n)$ on the coordinate plane, we obtain the graph of a function 
continuous on $(0,\infty)$ (Theorem \ref{ErgMeanStab}). The behavior of this function in 
the neighborhood of the origin is more complicated. We show in Example 3 below 
the existence of an $S$-integrable function, for which there exist very 
big intervals $\bar{m},\bar{n}$ such that $\frac n{|Y|}\approx\frac m{|Y|}\approx 0$, 
but $A_n\not\approx A_m$. However, Theorem \ref{NSBET} shows stabilization of ergodic means
on some initial segment of very big moments. In other words there exists a very big number 
$m$ such that for all very big numbers $n<m$ one has $A_n\approx A_m$ for almost all 
$y\in Y$, i.e. the share of those $y\in Y$, for which the statement is not true, is
infinitesimal. It is interesting that the proof of this theorem uses the Ergodic Theorem 
for infinite probability spaces and is equivalent to the last one in some sense.

We consider specially the case of discrete dynamical systems that are approximations of dynamical systems on compact metric spaces. We introduce here a definition of
such approximations (Definition \ref{hypap} below). The existence of approximations in the sense of Definition \ref{hypap} is proved
for a huge class of dynamical systems on compact metric spaces (see Section 4).

The approach to approximation suggested here differs from the most popular approach in ergodic theory based on Rokhlin's Theorem (see e.g. \cite{KSF}). The Nonstandard Analysis (NSA) approach to Rokhlin's finite approximations of Lebesgue dynamical systems will be discussed in another paper. Some preliminary results were announced in
\cite{GHL}. Rokhlin's approximations have many interesting applications to ergodic theory, especially to problems connected with the entropy of dynamical systems. However, Definition \ref{hypap} is more appropriate for investigation of computer simulation of continuous dynamical systems (see e.g. Example 6). We show also that the existence of a dynamical system on a compact metric space, for which a given very big finite dynamical system is an approximation, gives some additional information about the the behavior of the finite dynamical system on very big intervals of time (see Proposition \ref{functional} and Theorem \ref{un-erg} below).

Our approach provides some deeper understanding of the interrelation between very big discrete dynamical systems and continuous dynamical systems it the spirit of the approach formulated in \cite{Ze}: "Continuous analysis and geometry are just degenerate approximations to the discrete world... While discrete analysis is conceptually simpler ... than continuous analysis, technically it is usually much more difficult. Granted, real geometry and analysis were necessary simplifications to enable humans to make progress in science and mathematics....". In some sense, our paper contributes to this idea for dynamical systems.

Properties involved in the discussion above (very big set, very small number, etc.) obviously are not well defined.
They strongly depend on the problems, where they are used. Let us call them \emph{vague} properties. They cannot be formalized in the framework of the standard mathematics based on the G. Cantor's Set Theory. In this theory a set is understood as a collection of objects that satisfy a certain property that is well defined. This means that any two persons agree about any object, if this object has a given property or not. In other words one can definitely say about any object, if this object is an element of a given set or this is not the case.

Vague properties do not define sets. Consider, for example, the collection 
$\Omega$ of all very big natural numbers. If we accept the existence of a very big 
number, then the collection $\Omega$ is not a set. Indeed, if $\Omega$ is a set, then, 
obviously  $\N\setminus\Omega$ is a set. It is clear intuitively, that 
$0\in\N\setminus\Omega$ and if $n\in\N\setminus\Omega$, then  
$n+1\in\N\setminus\Omega$. So, by Axiom of Induction $\Omega=\emptyset$.

This argument comes up to the well-known paradox of a heap sand due to Eubulides (IV century B.C.): A single grain of sand is certainly not a heap. Nor is the addition of a single grain of sand enough to transform a non-heap into a heap: when we have a collection of grains of sand that is not a heap, then adding a single grain will not create a heap. And yet we know that at some point we will have a heap.

This paradox cannot be resolved in the framework of conventional (standard) mathematics, since the property "to be a heap of sand" is a vague property. On the other hand, vague properties are very common in natural sciences, economy and other areas of application of mathematics. Arguments, using them can be met in many investigations in these areas. These arguments seem to be quite convincible. Moreover, we will see below that the formalization in the framework of standard mathematics of some statements and arguments involving vague properties may be too complicated or even irrelevant.

Nonstandard Analysis, discovered by A. Robinson in the 60-s of the previous century introduced constant infinitesimals and infinitely large numbers in mathematics on the contemporary level of mathematical rigor. It opened the way to use vague collections (called the external sets in NSA) rigorously. The methods of NSA found numerous applications in the various areas of mathematics from mathematical physics to mathematical economics
(see, e.g. \cite{Alb, GKK, Wolf}). However, in the most of the papers nonstandard analysis is used as a tool to obtain results in standard mathematics.

The results mentioned in the beginning of this Introduction have more natural formulations in terms of vague properties rather than in the framework of the standard mathematics. Some of them, like Theorem \ref{ErgMeanStab}, can be simply reformulated in the framework of standard mathematics in terms of sequences of finite probability spaces, other, like Theorem \ref{NSBET}, do not have simple meaningful standard formulation. However, Theorem \ref{NSBET} has clear meaningful sense and can even be monitored in computer experiments (see Example 3).

Section 2 contains a brief introduction to NSA.  We discuss the formalization of the vague properties mentioned above. In particular, the formalization of the notion of a very big (very small) number is the formal definition of \emph{infinite (infinitesimal)} numbers in the NSA. The numbers that are not very big are called \emph{bounded} or \emph{finite}. We say that two elements $\a$ and $\b$ of a metric space are infinitesimally close ($\a\approx\b$), if the distance between them is  infinitesimal. A very big finite set is defined as a set, whose cardinality is an infinite natural number. Not very big sets are said to be standardly finite. As a rule we call them just finite sets if it does not yield a misunderstanding. There are no new results in this section, however, the exposition of the introduction to NSA is new. Some proofs in this section are given for illustration of the basic principles of NSA.

We try to make the exposition in Section 2 not too formal. Nonstandard analysis has one feature, that makes the achievement of this goal a little bit more difficult, than in the standard mathematics. One of the main axioms of NSA, the \emph{Transfer principle} is a statement about all conventional mathematical propositions. To make this statement mathematically rigorous one has to provide a formal definition of mathematical language. This can be done, for example, in the framework of the Axiomatic Set Theory. To avoid excessive formalization, we skip formal description of the mathematical language, assuming that every reader understand intuitively  what is a conventional mathematical proposition and can easily apply the Transfer Principle to any concrete proposition.

The formulations of the main results of the paper formulated both in the language of the NSA and in the framework of the standard mathematic, their discussion, illustration by examples and proofs of some simple statements are contained in Section 3. It is possible to understand the formulations of the main theorems of the paper (Theorem \ref{ErgMeanStab}, \ref{NSBET}, \ref{cycle-appr} and \ref{un-erg}) on the intuitive ("physical") level before reading 
Section 2,  if to interpret an infinite number as a very big number, an infinitesimal as a very small number, a hyperfinite set as very big finite set, an internal set as a usual set, an external set
as a collection of objects defined by a vague property. For example, under this interpretation the set $\*\N$ is the set of all natural numbers, while $\N$ is a collection of not very big numbers.

The rigorous formulations of the definitions and of the main results in terms of sequences 
of dynamical systems on finite probability spaces are contained in the part iv) 
of Section 3. The proofs of Theorems~\ref{ErgMeanStab},\ref{NSBET},\ref{cycle-appr} are 
contained in Section~4.

\textbf{Acknowledgements}. The authors are grateful to Peter Loeb, Edgardo Ugalde for 
helpful discussions of various parts of this paper and to Andrew Mertz and Kamlesh Parwani 
for their help with computer experiments.

\section{Basic Nonstandard Analysis}

i) We deal with some standard universe $\S$ that contains all objects necessary to develop a huge part of standard mathematics.

\begin{Def} \label{st-un}
A set $\S$ is said to be a standard universe, if
\begin{enumerate}
\item the field $\R\in\S$,
\item $a\in A\in\S\Limpl a\in\S$,
\item $A\in\S\Limpl\P(A)=\{B\ss A\}\in\S$,
\item $A,B\in\S\Limpl A\times B, A^B\in\S$,
\item if all elements $a$ of $A\in\S$ are sets, then $\left(\bigcup\limits_{a\in A}a\right)\in\S$,
\item any finite set of elements of $\S$ is an element of $\S$,
\end{enumerate}
\end{Def}

\begin{Prop} \label{st-prod} 1). If a set $A\in\S$ and $B\ss A$, then $B\in\S$

2). If $I\in\S$, $\{A_i\ |\ i\in I\}\in\S$, then $\prod\limits_{i\in I}A_i\in\S$\end{Prop}

\textbf{Proof}. 1) Since $B\in\P(A)$ and $\P(A)\in\S$ by property~3 of Definition \ref{st-un} then $B\in\S$ by property 2.

2) The set  $B=\prod\limits_{i\in I}A_i\in\S\ss\left(\bigcup\limits_{i\in I}A_i\right)^I=A$. The set $A\in\S$ by properties 4 and 5. Thus, $B\in\S$ by the statement 1 of this 
Proposition. $\Box$

We use the following notation: Let $\frak{A}$ be any collection of objects and let $\Phi$ be a standard sentence. We write $\frak{A}\models\Phi$ if $\Phi$ is true in
$\frak{A}$. This means that $\Phi$ is true, if all variables involved in it assume values in $\frak{A}$. In this paper we use for $\frak{A}$ either the collection $\S$ or the collection $\*\S$ introduced below.

It is easy to see that most part of the mathematical theorems are true in $\S$. Indeed, since the field $\R\in\S$ then the operations of addition and multiplication, as well as the order relation is in $\S$. Since the elements of all sets in $\S$ are also in $\S$ the axioms of linearly ordered fields for $\R$ are true in $\S$. The axiom
of the least upper bound for $\R$ is true in $\S$ due to Proposition \ref{st-prod} (1). One of the important axioms of set theory is the

\emph{Separation Axiom}. For an arbitrary set $B$, standard property $\Phi(x,y_1,...y_n)$ and elements $t_1,...t_n$ there exists a set $C=\{b\in B\ |\ \Phi(b,t_1,...,t_n)\mbox{is true}\}$.

This axiom also follows from Proposition \ref{st-prod} (1): if $B\in\S$ and $t_1,...t_n\in\S$, then $C\in\S$, since $C\ss B$.

The \emph{Axiom of Choice} states that the direct product of any family of non-empty sets is a non-empty set. Its validity in $\S$ follows from Proposition \ref{st-prod} (2) and Definition \ref{st-un} (1).

It is accepted by the most part of mathematicians that any mathematical statement can be formalized and proved (if it is provable) in the framework of the Zermelo-Fraenkel axiomatics for set theory ({\bf ZFC}). Besides the Separation Axiom and the Axiom of Choice this system contains the Axiom of infinity that is true in $\S$, since $\N\in\S$, the Axiom of the Unordered Pair that follows from Definition \ref{st-un} (6), the Axiom of Union, the Axiom of the Power Set that follows from Definition \ref{st-un} (3) , the Axiom of Regularity~\footnote{Actually the axiom of regularity was introduce later by John von Neumann} and the Axiom of Replacement, that is not true in $\S$. The axiomatics that contains all listed above axioms except the Axiom of Replacement is the Zermelo axiomatics. The Zermelo Axiomatic is enough for formalization of all concrete mathematics (analysis, differential equations, mathematical physics, geometry, etc.), while the Replacement Axiom is used only for the needs of Foundations of Mathematics. For example, it is used in the proof of existence of a set $\S$, that satisfies Definition \ref{st-un}.

The above discussion justifies the following

\textbf{Metatheorem 1}. Every theorem provable in Zermelo Axiomatics is true in $\S$.

\bigskip

ii) We extend the standard universe $\S$ by adding infinite, numbers, infinitesimals and 
some other objects. A good intuition for working with the nonstandard
extension $\*\S$ of $\S$ is provided by the following point of view. We consider the 
standard universe $\S$ as the universe of visual objects, while $\*\S$ is obtained by 
adding to $\S$ objects visual through a microscope (e.g. infinitesimals) and through 
a telescope (e.g. infinite numbers).

If $t\in\S$ and $t$ is not a set, then $t\in\*\S$. If a set $A\in\S$, then $A$ may be 
extended in $\*\S$, by adding some nonstandard elements. The
nonstandard extension of a set $A\in\S$ is denoted by $\*A$. The set $\*A\in\*\S$. 
For example, we will see later that the nonstandard extension $\*\R$ of the set $\R$ 
consists of infinite numbers, infinitesimals and numbers of the type $t+\a$, where $t\in\R$
and $\a\approx 0$. The elements of $\S$, that are not sets and the sets of the the form 
$\*A$, where $A\in\S$, are said to be \textbf{standard} elements of $\*\S$. To study the 
universe $\*\S$ we use the conventional language of mathematics extended by the predicate 
$\bs(x)$ that is interpreted as "$x$ is standard". We denote by $\bs$ also the collection 
of all standard elements of $\*\S$. So, to write $\bs(x)$ is the same as to write 
$x\in \bs$. We use the abbreviations $\alls\,x\ ...$ and $\exs\, x\ ...$ for 
$\all\,x (\bs(x)\Limpl ...)$ and $\ex\,x\ (\bs(x)\& ...)$ respectively. Let $\Phi$ be 
a proposition that may contain some free variables assuming values in $\S$ or in $\*\S$. 
Then $\Phi\st$ is a proposition that is obtained from $\Phi$ by replacing any quantifier 
$\all\ (\ex)$ by $\all\st\ (\ex\st)$. All elements of $\*\S$ (sets and not sets) are said 
to be internal elements. Propositions formulated in conventional language are said to be 
\emph{internal}. Propositions containing the predicate $\bs$ and the map $\*$ are said 
to be \emph{external}. External propositions are used to describe vague properties 
discussed in the Introduction.

We introduce now the axioms for the nonstandard universe $\*\S$. These axioms are simplified versions of the axioms of one the axiomatic nonstandard set theories - the theory \textbf{HST} (Hrbacek Set Theory) (\cite{KR}).

We say that a proposition is a \emph{sentence} if all variables involved in $\Phi$ are connected by quantifiers

I. There exist an injective map $\*:\S\to\bs$ ($\*(\S):=\bs$) such that for any element $t\in\S$ that is not a set one has $\*t=t$.

II. (Transitivity of $\*\S$) If $A\in\*\S$ and $a\in A$, then $a\in\*\S$.

III. (Transfer Principle) If $\Phi(x_1,...,x_n)$ is an internal proposition and $a_1,...,a_n\in\S$, then
$$\S\models\Phi(a_1,...,a_n)\Liff\*\S\models\Phi(\*a_1,...,\*a_n)\Liff\*\S\models\Phi\st(\*a_1,...,\*\a_n).$$

The Transfer Principle immediately implies

\textbf{Metatheorem 2}. Every theorem provable in Zermelo Axiomatics is true in $\*\S$.

The next propositions easily follow from the Transfer Principle.

\begin{Prop} \label{transfer1}  The bijection $\*$ preserves the boolean operations on sets and finite cartesian products. \end{Prop}
\textbf{Proof}. Let $A,B,t\in\S$. Then one has $t\in(A\cap B)\liff t\in A\& t\in B$. So, by Transfer Principle
$$
\*t\in\*(A\cap B)\liff \*t\in\* A\& \*t\in \*B.
$$
Thus, 
$$
\all\st t(t\in\*(A\cap B)\liff t\in\*A\& t\in\*B).
$$ 
Again, by Transfer Principle $$\all t(t\in\*(A\cap B)\liff t\in\*A\& t\in\*B)$$
is true, which means that $\*(A\cap B)=\*A\cap\*B$. For the other operations the proof is similar. $\Box$

\begin{rem} To understand the second part of the proof, keep in mind that the set $\*A$, $\*B$ and $\*(A\cap B)$ may contain not only standard elements. \end{rem}

We say that the proposition $\Phi(x)$ defines the standard element $t\in\S$ if the following statement is true in $\S$:
$$\Phi(t)\&\all\, y(\Phi(y)\to y=t)$$.

\begin{Prop} \label{transfer2} If a proposition $\Phi(x)$ defines an element $t\in\S$, then it defines $\*t\in\*\S$\end{Prop} $\Box$
\begin{Cor} \label{t2} a). $\*\emptyset=\emptyset$.

b). If a set $A=\{\a_1,...,a_n\}\in S$, then $\*A=\{\*a_1,...,\*a_n\}$.
\end{Cor}

To prove the statement b) of this Corollary first prove it for $n=1$ using the Transfer Principle, then apply the induction by $n\in\N$.

IV. (Idealization Principle) If a set $A\in\S$ is infinite, then $\*A\setminus A\neq\emptyset$.

\begin{Prop} \label{infinite} If $N\in\*\N\setminus\N=\iN$, then for any $n\in\N$ one has $N>n$.\end{Prop}
\textbf{Proof}. If $N\leq n$ for some $n\in\N$, then $N\in\*\{0,....n\}=\{0,...,n\}$ 
by Corollary \ref{t2} (b). Thus, $N\in\N$. The contradiction. $\Box$

Obviously, if $N\in\iN$, then $N-1\in\iN$. Thus, the set $\iN$ does not have a minimal element and the set $\N$ satisfies the antecedent of Induction Principle, however
$\N\neq\*\N$. So, the sets $\N$ and $\iN$ are not internal sets, since, by the Transfer Principle, the Induction Principle is applicable to internal subsets of $\*\N$.

We see, thus, that the property 3 of Definition \ref{st-un} fails for the nonstandard universe $\*\S$. There exists a set $A$ and a subset $B\ss A$, such that
$A\in\*\S$ and $B\notin\*\S$.

\begin{Def} \label{external} We say that a set $B$ is external, if it is not an internal set, but is a subset of an internal set. \end{Def}

We extend the nonstandard universe by adding all external sets: 
$\ws:=\*\S\cup\{B\ss\*S\ |\ B\ \mbox{is an external set}\}$

We use the abbreviations $\all\il\,x\ ...$ and $\ex\il\, x\ ...$ for $\all\,x (x\in\*\S\Limpl ...)$ and $\ex\,x\ x\in\*\S\& ...)$ respectively. Let $\Phi$ be a proposition (maybe external) that may contain some free variables assuming values in $\S$ or in $\*\S$. Then $\Phi\il$ is a proposition that is obtained from $\Phi$ by replacing any quantifier $\all\ (\ex)$ by $\all\il\ (\ex\il)$. There is no need to write ${\all\il}\st x$ or ${\ex\il}\st x$, since if $x\in\bs\Limpl x\in\*\S$.  
Notice that $\Phi\il$ is an \emph{external} proposition even if$\Phi$ is an internal proposition since the proposition $x\in\*\S$ is defined in terms of 
the map $\*$. Obviously for every external proposition $\Phi$ and every internal 
proposition $\Psi$ one has
\begin{equation}\label{notations}
\ws\models\Phi\il\Liff\*\S\models\Phi,\quad\ws\models\Psi\st\Liff\S\models\Psi.
\end{equation}

\begin{rem}\label{nelson} The propositions of the form $\Phi\il$ we call in this paper \emph{Nelson-type propositions}, since E. Nelson was the first who suggested a formal axiomatic (IST) for internal sets in conventional language extended by the predicate "$x$ is standard" - the analog of our predicate $x\in\bs$ and who wrote the first exposition of probability theory in the framework of IST \cite{N1,N2}

The main results of this paper are formulated as Nelson-type sentences, since only this sentences are most intuitively clear formalizations of the statements containing the vague notions discussed above. However, in the proves we use propositions that involve variables assuming values of arbitrary external sets. These propositions make
many of proofs much simpler, than if we restrict ourselves to the arguments that can be formalized in IST.
\end{rem}

\begin{Def} \label{finite} We say that an internal set $A$ is finite, if there exists $n\in\*\N$ and an internal bijection $\f:A\to\{1,...,n\}$. In this case we say that the cardinality $|A|$ of $A$ is equal to $n$. If $|A|\in\iN$, then we say that $A$ is a hyperfinite set. If $|A|\in\N$, then we say that $A$ is standardly finite (s-finite) set. \end{Def}

\begin{rem} \label {st-card} If $\Phi$ is the definition of a finite set in the conventional mathematics, then the statement used for the definition of a finite set in Definition \ref{finite} is the proposition $\Phi\il$. So it would be more correct to call the number $|A|$ defined in Definition \ref{finite} the internal cardinality.
However, the real cardinality of $A$, i.e. the cardinality of $A$ in the "global" universe of all sets may strongly depend on the properties of $\*\S$ and is never used in applications of NSA. So, we prefer to to keep the term "cardinality" for the internal cardinality, and to call the real cardinality of a set the external cardinality. It agrees with our intuition,
according to which the set $\*\N$ is the set of all natural numbers that includes also 
those numbers that can be seen through a telescope.
The cardinality should be a well-defined (not vague) notion, that is why it must be 
an internal (not external sets). The definition of a hyperfinite set is a formalization of the vague notion of a very big set. The definition of an s-finite set is a formalization of the vague notion of a not too big set. Obviously the external cardinality of a hyperfinite set is infinite. On the other hand, it can be easily proved by induction that if an internal set is s-finite, then its internal cardinality is equal to its external cardinality and every set, whose external cardinality is a standard natural numbers, is an internals set.\end{rem}

V. (Saturation Principle)~\footnote{We introduce here the weakest form of the Saturation  Principle. However, this form is enough for our goals.}. If an external sequence $\{A_n\ |\ n\in\N\}$ of internal sets has a finite intersection property (i.e. $\all n\in\N\ \bigcap\limits_{k\leq n}A_k\neq\emptyset$), then $\bigcap\limits_{n\in\iN }A_k\neq\emptyset$.

\bigskip

iii) The following definition contains the formalization of the vague notions of very big and very small numbers.

\begin{Def}\label{infin}  We say that
\begin{enumerate}
\item a number $\Omega\in\*\R$ is infinite ($\Omega\sim\infty$), if $|\Omega|>N$ for all 
$N\in\N$. A number $\a\in\*\R$ that is non-infinite is said to be bounded or finite 
($\a\ll\infty$).
\item A number $\a\in\*\R$ is said to be infinitesimal $\a\approx 0$, if $|\a|<\frac 1N$ 
for all $N\in\N$. Two numbers $\a$ and $\b$ are infinitesimally close $\a\approx\b$, if $\a-\b\approx 0$. We write $|a|\gg 0$, if $\a$ is not an infinitesimal number.
\item A number $t\in\R$ is said to be a standard part (or a shadow) of a bounded number 
$\a$ ($t=\o\a$), if $t\approx \a$.
\end{enumerate}
\end{Def}

The existence of infinite and infinitesimal numbers follows from Proposition \ref{infinite}

We denote the set of all bounded numbers by $\*\R_b$. If $t\in\R$, then the set 
$\Mu(t)=\{\a\in\*\R\ |\a\approx t\}$ is called the \emph{monad} of $t$.

The properties of infinite, bounded and infinitesimal numbers are similar to the properties of sequences that diverge to infinity, are bounded and tend to $0$ respectively in standard calculus. They can be summarized as follows.

\begin{Prop} \label{inf-prop} 1)  $\Omega\sim\infty\Liff \Omega^{-1}\approx 0$.

2) The set $\*\R_b$ is a subring of the field $\*\R$ and the set $\Mu(0)$ is an ideal in the ring $\*\R_b$.
\end{Prop}

We leave a simple proof of this proposition as an exercise.

\begin{Th}\label{st-part} 1) Every $\a\in\*\R_b$ has a unique standard part.

2) The map $\o:\*\R_b\to\R$ is a homomorphism of a ring $\*\R_b$ onto the field $\R$.
\end{Th}

\textbf{Proof.} The only non-trivial statement is the existence of standard part for any 
bounded element.  Let $\a\in\*R_b$ be bounded. Then there exists $s\in\R$ such that $s>\a$.
Consider the set $P=\{p\in\R\ |\ p<\a\}\ss\R$. This set is nonempty since $\a$ is bounded not only from above, but also from below. Thus, there exists
$t=\sup P$ in $\R$. If $t<\a$, then $t+\frac 1n<\a$ for any $n\in\N$. Otherwise $\sup P\geq t+\frac 1n$ for some $n\in\N$. So, $\a-t<\frac 1n$ for any $n\in\N$, i.e. $\a-t\approx 0$. If $t>\a$ the proof is similar. $\Box$ \renewcommand{\sp}{\mbox{st}}

Let $(X,\rho)\in\S$ be a metric space. In what follows we write this and similar sentences as "Let $(X,\rho)$ be a standard metric space". Then by the Transfer Principle $\*\S\models (\*X,\*\rho)\ \mbox{is a metric space}$. In what follows for any proposition $\Phi$
instead of writing $\*\S\models\Phi$ we write "$\Phi$ in $\*\S$. For example, the previous statement may be written as $(\*X\*\rho)$ is a metric space in $\*\S$.

For any $\xi_1,\xi_2\in\*X$ we write $\xi_1\approx \xi_2$, if $\*\rho(\xi_1,\xi_2)\approx 0$. For $x\in X$ and $\xi\in\*X$ we write $(x=\sp(\xi)$, if $x\approx\xi$. We say in this case that $x$ is a standard part of $\xi$~\footnote{For the case of the metric space $\R$ we use also the notation $\o$ for the standard part.}. Obviously, $\sp(\xi)$ is defined uniquely. An element of $\xi\in\*X$ is said to be \emph{nearstandard}, $\sp(\xi)$ exists. 
In particular, $\*\R_b$ is the set of all nearstandard elements of $\*\R$. Similarly, the set of all nearstandard elements of $\*X$ is denoted by $\*X_b$. For $x\in X$ the set $\{\xi\in\*X\ |\ \xi\approx x\}$ is said to be a monad of $x$ and denoted by $\Mu(x)$. For an arbitrary $0<\e$ and $x\in X$ let $B_{\e}(x)=\{\xi\in X\ |\ \rho(\xi,x)\leq\e\}$. Then it is easy to see that
\begin{equation}\label{monad}
\Mu(x)=\bigcap\limits_{n\in\N}\*B_{\frac 1n}(x).
\end{equation}

The simple proof of the following proposition can be found in the books \cite{Alb, Wolf, GKK}.

\begin{Prop} \label{metric} 
Let $(X,\rho)$ be a standard separable metric space and $A\ss X$. 
Then the following statements are true~\footnote{These statements are true without the 
assumption of separability of $X$, if the nonstandard universe $\*X$ satisfies some 
stronger Saturation Principle.}.
\begin{enumerate}
\item The set $A$ is open if and only if $\all\,x\in A\ \Mu(x)\ss\*A$.
\item The set $A$ is closed if and only if $\all\, x\in X\ (\Mu(x)\cap\*A\neq\emptyset\Limpl x\in A)$.
\item The set $A$ is compact if and only if it is closed and every element $\xi\in\*A$ is nearstandard. In particular, $X$ is a compact metric space, if and only if every element $\xi\in\*X$ is nearstandard .
\end{enumerate}
\end{Prop}
\textbf{Proof}. To illustrate how the axioms I-V work let us prove that every $X$ compact metric space satisfies the second statement of 3). Suppose that $X$ is a compact metric space. The for every $n\in\N$ there exists a finite set $\{x_1,...,x_k\}\ss X$ such that $X=\bigcup\limits_{i=1}^k B_{\frac 1n}(x_i)$. Then, by the Transfer Principle 
$\*X=\bigcup\limits_{i=1}^k \*B_{\frac 1n}(x_i)$. In other words, 
$\all\,\xi\in\*X\,\all\, n\in\N\,\ex\st x\ \xi\in\*B_{\frac 1n}(x)$ (A). 
Suppose that there exists $\xi\in\*X$ that is not nearstandard. It means that 
$\all\st x\ \xi\notin\Mu(x)$, i.e. by formula (\ref{monad}) 
$\xi\notin\bigcap\limits_{n\in\N}\*B_{\frac 1n}(x)$. Then, by the Saturation Principle 
there exists $n\in\iN$ such that $\xi\notin\*B_{\frac 1n}(x)$, which contradicts (A). $\Box$

iii) This part of Section 2 contains some well-known facts of nonstandard analysis that are necessary only for the proofs of results formulated in Section 3. These proofs are contained in Section 4.

\begin{Th}\label{spill} Let $A\ss\*\R$ be an internal set.
\begin{enumerate}
\item If $\N\ss\A$, then $\{0,1,...,N\}\ss A$ for some $N\in\iN$.
\item If $\iN\ss\A$, then $\*\N\setminus\{0,1,...,n\}\ss A$ for some $N\in\N$.
\item If $\Mu(0)\ss A$, then $(-t,t)\ss A$ for some $0<t\in\R$.
\item (Robinson's Lemma) Let $\la s_n\ |\ n\in\*\N\ra$ be an internal sequence such that $s_n\approx 0$ for all $n\in\N$. Then there exists $N\in\iN$ such that
$s_n\approx 0$ for all $n<N$
\item Let $\la N_n\ |\ n\in\N\ra$ be an external sequence of infinite numbers. Then there exists $N\in\iN$ such that $N<N_n$ for all $n\in\N$.
\end{enumerate}
\end{Th}

Proofs of the statements 1-4 can be found e.g. in \cite{Wolf}, p. 53.

\textbf{Proof of the statement 5.} Consider an external sequence of internal sets 
$B_n\la\{n,...,N_n\}\ |\ n\in\N\ra$. This sequence has obviously a finite intersection 
property. By the Saturation Principle $\bigcap\limits_{n\in\N}B_n\neq\emptyset$. 
Any element of this intersection satisfies the conditions of the statement 5. $\Box$

We say that a set $I\ss\iN$ is an \emph{initial segment of infinite numbers} if $I=\{0,...,N\}\cap\iN$.

\begin{Th} \label{sequences}   
Let $\la a_n\ |\ n\in\*\N\ra$ be an internal sequence of nonstandard real numbers 
(elements of $\*\R$). Then $\lim\limits_{n\to\infty}\o a_n=a\in\R$ if and only if 
$a_L\approx a$ for all $L$ in some initial segment of infinite numbers. 
\end{Th}

\textbf{Proof}. 
$Longrightarrow$ Let $\lim\limits_{n\to\infty}\o a_n=a\in\R$. For any $m\in\N$ there 
exists $n_m\in\N$ such that the internal set $A_m=\{k\in\*\N\ |\ |a_k-a|\leq\frac 1m\}$ 
contains all $k\in\N$ such that $k\geq n_m$. Then, by Theorem \ref{spill} (1) there exist 
$N_m\in\iN$ such that $\{n_m,...,N_m\}\subset A_m$. By Theorem \ref{spill} (1) there 
exists $N\in\iN$ such that $N<N_m$ for all $m\in\N$. This $N$ satisfies the conditions of 
the theorem. 

$\Longleftarrow$ Let $N\in\iN$ be such that $a_L\approx a$ for all $L\leq N,\ L\in\iN$. 
Fix an arbitrary $m\in\N$. Then the internal set
$B=\{k\in\*\N\ |\ |a_k-a|<\frac 1m\}\supseteq\{L\in\iN\ |\ L\leq N\}$. Thus, by Theorem \ref{spill} (2) there exists $n\in\N$ such that
$B\supseteq\{n,...,N\}\supset\{k\in\N\ |\ k>n\}$. So, $\all\, m\in\N\,\ex\, n\in\N\,\all\,k\in\N\ k>n\Limpl |\o a_k-a|<\frac 1m$. This means that $\lim\limits_{k\to\infty}=a$. $\Box$

\begin{rem} \label{st-nonst} Notice, that the sequence $\la \o a_n\ |\ n\in\ N\ra\in\S$ 
is a standard sequence and, thus, the sequence $\la\*\o a_n\ |\ n\in\*\N\ra$ is defined. 
This sequence \emph{is not necessarily equal} to the initial internal sequence  
$\la a_n\ |\ n\in\*N\ra$.  The only statement that can be claimed is that the entrees 
of these two sequences are infinitesimally close on an interval $\{0,...,N\}$ for some 
$N\in\iN$.
\end{rem}

Let $X$ be a compact metric space, $Y\ss\*X$. We say that $Y$ is a \emph{dense} subset of $\*X$, if  $\all\, x\in X\ex\,y\in Y\ x\approx y$. Proposition \ref{metric} implies that the last statement is equivalent to the statement $\all\, x\in\*X\,\ex\,y\in Y\ x\approx y$.

If $Y\ss X$ and $X$ is a metric space, then we say that a function $f:Y\to\*\R$ is 
$S$-continuous on $Y$, if $\all\, y_1,y_2\in Y\ y_1\approx y_2\Limpl f(y_1)\approx f(y_2)$.

\begin{Th} \label{S-cont} Let $(X, \rho)$, $(Z, d)$ be standard separable metric spaces, 
$X$ be a compact space and $Y\ss\*X$ be a dense internal subset of $\*X$.

1) A function $f:X\to Z$ is continuous if and only if $\*f\*X to\*Z$ is $S$-continuous on 
$\*X$.

2) Let $F:Y\to\*Z_b$ be an internal function that is $S$-continuous on $Y$, then the function $f:X\to Z$ defined by the formula $f(\sp(\xi))=\sp(F(\xi))$ is a continuous function.
\end{Th}

In what follows the function $f$ defined in the statement 2) of the theorem is said 
to be the \emph{visual image} of $F$, if $F\ss\*\R\times\*\R$. More generally,
if $A\ss\*\R\times\*R$, the the set $\o A=\{(\o a,\o b)\ |\ (a,b)\in A\}$ is said 
to be the visual image of $A$. This definition is specific for this paper. Usually, in NSA  
the set $\o A$ is said to be the shadow of $A$

\textbf{Proof}. We prove the statement $\Limpl$ for 1) and the statement $\Longleftarrow$ for 2).

1) $\Limpl$. Since $X$ is a compact space and $f$ is continuous on $X$, then $f$ is uniformly continuous on $X$. This means that
$\all\,\e\,\ex\,\d,\all\,x_1, x_2\in X\ \rho(x_1,x_2)<\d\Limpl d(f(x_1),f(x_2))<\e$ is true in $\S$. By the Transfer Principle the statement
\begin{equation} \label{S-cont1}
\all\st\,\e\,\ex\st\,\d,\all\,\xi_1, \xi_2\in\*X\ \*\rho(\xi_1,\xi_2)<\d\Limpl \*d(f(\xi_1),f(\xi_2))<\e
\end{equation}
is true in $\*\S$. If $\xi_1\approx\xi_2$, then the antecedent of the implication in the statement (\ref{S-cont1}) is true for any standard $\d$. Thus, the consequence of this implication is true for any standard $\e>0$. This means that $\*f(\xi_1)\approx\*f(\xi_2)$.

2) $\Longleftarrow$. Due to Proposition \ref{metric} (3) the function $f$ is defined 
correctly. We have to prove that $f$ is uniformly continuous on $X$.
By the definition of $f$ it is enough to prove that
\begin{equation} \label{S-cont2}
\all\st\e>0\,\ex\st\d>0\,\all\, \xi_1,\xi_2\in Y\ \*\rho(\xi_1,\xi_2)<\d\Limpl \*d(F(\xi_1),F(\xi_2))<\e.
\end{equation}

Fix an arbitrary standard $\e>0$. Due to the $S$-continuity of $F$, the internal set
$$B=\{0<\d\in\*\R\ |\ \all\,\xi_1,\xi_2\in Y\ \*\rho(\xi_1,\xi_2)<\d\Limpl \*d(F(\xi_1),F(\xi_2))<\e\}$$
contains all $0<\d\approx 0$. Thus, by Theorem \ref{spill} (3), there exist a $\d\in\S$ such that $\*(0,\d)\ss B$. This proves the statement (\ref{S-cont2}). $\Box$

\bigskip

iv) We list now the necessary definitions and facts concerning Loeb spaces.
We need here only a particular case of a Loeb space, namely the Loeb space constructed from the
hyperfinite set $Y$ endowed with the uniform probability measure.

Define the internal finitely additive measure $\mu$ on the algebra $\P\il(Y)$ of internal subsets
of $Y$ by the formula
$$
\mu(B)=\frac{|B|}M.
$$

This measure induces the external finite additive measure $\o\mu$ on $\P\il(Y)$

The Saturation Principle and the Caratheodory Theorem imply the
possibility to extend $\o\mu$ on the $\s$-algebra $\s(\P\il(Y))$
generated by $\P\il(Y)$. The Loeb space with the underlying set
$Y$ is the probability, space $(Y,P_L(Y),\mu_L)$, where $P_L(Y)$
is the completion of $\s(\P\il(Y))$ with respect to the extension of $\o\mu$ and
$\mu_L$ is the extension of $\o\mu$ on $P_L(Y)$. The measure $\mu_L$
is said to be the \emph{Loeb measure} on $Y$. If necessary we use the notation
$\mu_L^Y$.  We need the following property of the Loeb measure that follows immediately 
from the Saturation Principle.

\begin{Prop} \label{Loeb-mes} For every set $A\in\P_L$ there exists an internal set 
$B\ss Y$ such that $\mu_L(A\Delta B)=0$. 
\end{Prop}

\begin{Cor}\label{a.-e.} If $A\in\P_L$, then
$$\mu_L(A)=1\Liff \all^{st}\e>0\,\ex\,B\in\P^{int}(Y)\ (B\ss A\land\mu(B)>1-\e)$$
\end{Cor}

\begin{rem} The proposition in the right hand side of this corollary is a Nelson-type 
proposition that will be used as a formalization of the notion "almost
everywhere in $Y$" ("for almost all $y\in Y$").
\end{rem}

For an arbitrary complete separable metric space $R$ and an external
function $f:Y\to R$ an internal function $F:Y\to\*R$ is said to be
\emph{a lifting} of $f$ if $\mu_L(\{y\in Y \ |\ F(y)\approx
f(y)\})=1.$

\begin{Prop} \label{mes-func} A function $f:Y\to R$ is measurable
iff it has a lifting.
\end{Prop}

An internal function $F:Y\to\*\R$ is said to be
\emph{$S$-integrable} if for all $K\in\iN$ one has
\begin{equation}\label{S-int-1}
\frac 1M\sum\limits_{\{y\in Y\ |\ |F(y)|>K\}}|F(y)|\approx 0.
\end{equation}

We need the following properties of $S$-integrable functions.

\begin{Prop} \label{S-int}
1) An $S$-integrable function is almost everywhere bounded.

2) An internal function $F:Y\to\*\R$ is $S$-integrable iff
$Av(|F|)=\frac 1{M}\sum\limits_{y\in Y}|F(y)|$ is bounded and 
$\frac 1{M}\sum\limits_{y\in A}|F(y)|\approx 0$ for every internal $A\ss Y$ such that
$\frac{|A|}M\approx 0$.

3) An external function $f:Y\to\R$ is integrable w.r.t. the Loeb measure $\mu_L$
iff it has an $S$-integrable lifting $F$, in which case
$$
\int\limits_Yfd\mu_L=\o Av(F).
$$
\end{Prop}

We address readers to
\cite{Alb,Wolf} for the proofs of Propositions \ref{Loeb-mes}, \ref{mes-func} and \ref{S-int}.

\bigskip

\section{Formulation and Discussion of results}.

i). In the sake of convenience of the references we recall the formulation of classical 
G. Birkhoff Ergodic Theorem.(see e.g \cite{KSF,Br}).

\begin{Th} \label{BET} Let $(X,\Sigma, \nu)$ be a probability space and $T:X\to
X$ a measure preserving transformation and $f\in L_1(X)$. Denote
by
$$
A_k(f,T,x)=\frac 1k\sum\limits_{i=0}^{k-1}f(T^ix).
$$

Then

\begin{itemize}

\item[1)] there exists the function $\hat f(x)\in L_1(X)$ such
that $A_k(f,T,x)\to\hat f(x)$ as $k\to\infty$ a.e.;

\item[2)] the function $\hat f$ is $T$-invariant, i.e. $\hat
f(Tx)=\hat f(x)$ for almost all $x\in X$;

\item[3)] $\int_Xfd\nu=\int_X\hat fd\nu$.
\end{itemize}
\end{Th}

If $Y$ is a finite set, $|Y|=M$, then every function on $Y$ is integrable, whatever $\Sigma$ and $\nu$ are.
We restrict ourselves to the case of the uniform measure:
$\mu(A)=\frac{|A|}M$ for any set $A\ss Y$. Then any measure preserving transformation $T:Y\to Y$ is a bijection and
the integral of a function $F:Y\to\*\R$ is the average 
$Av(F)=\frac 1M\sum\limits_{y\in Y}F(y)$ of $F$. Theorem \ref{BET} is proved very easily 
in this case. We reformulate it as a statement in the nonstandard universe $\*\S$ assuming 
that $F$ is an internal function, $Y$ is a hyperfinite set and, thus, $M\sim\infty$.

For any $y\in Y$ denote the $T$-orbit of $y$ by  $\mbox{Orb}(y)$ and the period of $y$ by $p(y)$.
\begin{Prop} \label{INTBET}

For any $y\in Y$, if $n\gg p(y)$, then $A_n(F,T,y)\approx \hat F(y)$, where

$$\hat F(y)=\frac 1{|\mbox{Orb}(y)|}\sum\limits_{z\in\mbox{Orb}(y)}F(y)$$.
\end{Prop}

\begin{Cor} \label{CORINTBET}
If $n\gg M$, and $T$ is a cycle of length $M$, then $$\forall\, y\in Y\, \hat f(f)\approx \frac 1M\sum\limits_{y\in X}f(y)= Av(f)$$.
\end{Cor}

We leave a simple proof of this proposition as an exercise (see also the proof of Theorem \ref{ErgMeanStab} below).

ii). For the case of $M\sim\infty$ it is interesting to study the behavior of ergodic means for $n\sim\infty$ but such that
$\frac nM\ll\infty$. We start with simple examples.

\textbf{Example 1.} Let $Y=\{0,\dots,M-1\}$, $T:Y\to Y$ is defined by the formula $T(y)= y+1(\mod M)$ for any $y\in Y$.
Consider the function $F:Y\to\*\R$ such that
\begin{equation} \label{ex01}
F(y)=\left\{\begin{array}{ll}M,\ &\mbox{if}\ y\ \mbox{is even},\\ -M,\ &\mbox{if}\ y\ \mbox{is odd}.\end{array}\right.
\end{equation}

Then for any $y\in Y$ one has
\begin{equation} \label{ex1}
A_n(F,T,y)=\frac Mn\cdot\left\{\begin{array}{rll}0,\ &\mbox{if}\ n\ \mbox{is even},\ &y\ \mbox{is any},\\  1,\ &\mbox{if}\ n\ \mbox{is odd}, &y\ \mbox{is even}\\ -1,\ &\mbox{if}\ n\ \mbox{is odd}, &y\ \mbox{is odd}.                 \end{array}\right.
\end{equation}

\begin{figure}
\centering
\includegraphics*[viewport= 50 400 800 700, scale=.5]{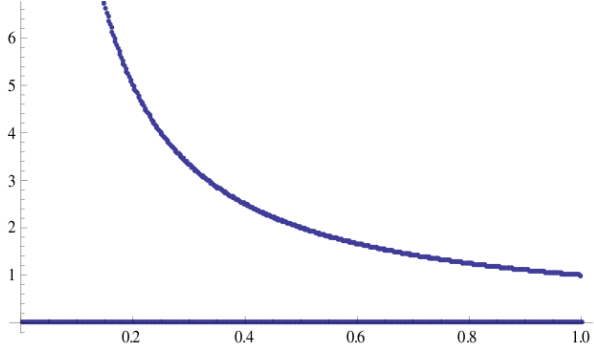}
\caption{$\Gamma(F)$ for $F$ defined in Eq.\ref{ex01}}
\label{fig1}
\end{figure}

Let us plot the set of points $\{(n/M, A_n(F,T,x))\ |\ n=0,...,kM\}$ for any chosen randomly $y\in Y$. The first question here is how to choose an infinite number $M$.
Recall that the notion of an infinite number is a formalization of a notion of a very big number. Certainly, the property "to be very big" depends on the problem. Sometimes very moderate numbers can be considered as very big. In this example we consider a number $M$ to be very big, if we see the the set of points $\{n/M\ | n=0,...,M\}$ as a continuous segment. 
On Fig.\ref{fig1} $M=1000$ and the randomly
chosen $y=698$.

We see on this picture the graphs of two functions $y=0$ and $y=\frac 1x$. The first one is the visual image of
of the set of points $A=\{(n/M, A_n(F,T,x))\ |\ n<M,\ n\ \mbox{is even}\}$, the second one is the visual image of
of the set of points $B=\{(n/M, A_n(F,T,x))\ |\ n<M,\ n\ \mbox{is odd}\}$.

Indeed, the set $A$ is an internal function, $\dom (A)=D_1=\{n/M\ |\ |\ n<M,\ n\ \mbox{is even}\}$ and $A(n/M)=0$. The set $D_1$ is dense in [0,1] and the function $A$ is obviously $S$-continuous on $D_1$. By Theorem \ref{S-cont}, the function $A$ defines the continuous function $f$. In this case obviously $f(t)\equiv 0$.

The set $B$ is an internal function, whose domain $D_2=\{n/M\ |\ |\ n<M,\ n\ \mbox{is odd}\}$ and $B(n/M)=\frac 1{n/M}$. The set $D_2$ is also dense in $\*[0,1]$, but
$B$ does not satisfy conditions of Theorem \ref{S-cont}, since $B(n/M)\sim\infty$ as $n/m\approx 0$. However, for any $0<a\in\R$ the function $B$ restricted to the set
$D_2\cap\*[a,1]$ is $S$-continuous. Obviously, its visual image is the function $f(t)=1/t$ restricted to the interval $[a,1]$. On Fig.\ref{fig1} $a=0.15$

In what follows we denote the set $\{(n/M, A_n(F,T,y))\ |\ 0<n/M\leq kM\}$ by $\Gamma_k(F)$. We write $\Gamma(F)$, if
$k=1$. \label{gamma}

It is natural to ask oneself the following question.

Under what conditions on an internal function $F:Y\to\*\R$ the visual image of the set $\Gamma_k(F)$ is the graph of a continuous function for any $k\ll\infty$ and for almost all $y\in Y$?

In view of the above discussion this question can be reformulated as follows.

\begin{Ques} \label{Q1}
Under what conditions on a function $F:Y\to\*\R$ the ergodic means $A_n(F,T,y)$ satisfy the following property:
\begin{equation}\label{cont-graph}
\frac nM\approx \frac mM\Limpl A_n(F,T,y)\approx A_m(F,T,y)
\end{equation}
for almost all  $y\in Y$?
\end{Ques}

In investigation of this question we restrict ourselves to the case, when a permutation $T:X\to X$ is a cycle of the length $M$. The general case can be easily reduced to this one.

Due to Corollary \ref{CORINTBET}, the implication (\ref{cont-graph}) holds for every function $F:Y\to \*R$ and every
$y\in Y$ if $\frac nM,\frac mM\sim\infty$. So, it is enough to consider the case of $\frac nM,\frac mM\ll\infty$.

In the following computer experiments, we illustrate that a proposition $\Phi$ holds for almost all $y\in Y$, by checking that this property holds for a randomly chosen
$y\in Y$, using computer generator of random elements. As in Example 1 we use a concrete very big finite set $Y$, that can be considered as hyperfinite one in our problem.

The property of $S$-integrability of a function on a very big finite space is an analog of the property of integrability of functions on infinite probability spaces. It is easy to see that any bounded function $F$ ($\max\limits_{y\in Y}|F(y)|\ll\infty$) is $S$-integrable. The $\d$-function gives an example of a function $F$ with bounded $Av(|F|)$ that is not $S$-integrable.

\textbf{Example 2} For the same $Y$ and $T$ as in Example 1 consider the function 
$F:Y\to\*\R$ given by the formula
\begin{equation} \label{delta}
F(k)=\left\{\begin{array}{ll} M,\ &k=0\ \\ 0,\ k\neq 0
\end{array}\right.
\end{equation}

\begin{figure} 
\begin{center}
\includegraphics*[viewport= 90 550 400 750, scale=0.5]{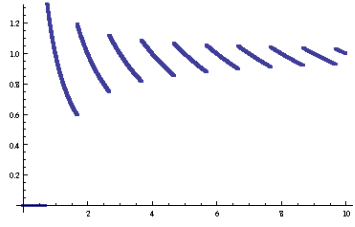}
\end{center}
\caption{$\Gamma_{10}(f)$}
\label{fig2}
\end{figure}

We leave to the reader as an easy exercise to find the formula for $A_n(F,T,y)$ for this function $F$. On Fig.\ref{fig2} we show the visual image of the sets
$\Gamma_{10}(f)$ for $M=1000$, and randomly chosen $x=322$. We see that the visual image on the first picture of Fig.\ref{fig2} is a graph of
a function that has points of discontinuity.

\begin{Th} \label{ErgMeanStab} Let $T:Y\to Y$ be a cycle of length $M$. Then
for every $S$-integrable function $F:Y\to\*\R$, for every
positive $a\in\R$ such that $0\ll a\ll\infty$ and for every numbers $K,L\sim\infty$
such that $\frac KM\approx\frac LM\approx a$ one
has $A_K(F,T,y)\approx A_L(F,T,y)$ for all $y\in Y$.
\end{Th}

\textbf{Proof}. Assume $K>L$ and estimate
$|A_K(F,T,y)-A_L(F,T,y)|$. It is easy to see that
$$
|A_K(F,T,y)-A_L(F,T,y)|\leq\left(\frac 1L-\frac
1K\right)\sum\limits_{k=0}^{L-1}|F(T^ky)|+\frac
1K\sum\limits_{k=L}^{K-1}|F(T^ky)|=U+V.
$$

One has
$$
U=\left(\frac ML-\frac MK\right)\frac
1M\sum\limits_{k=0}^{L-1}|F(T^ky)|\approx 0,
$$
since $\frac ML\approx\frac MK\approx\frac 1a$ and $\frac
1M\sum\limits_{k=0}^{L-1}|F(T^ky)|\leq\frac
1M\sum\limits_{k=0}^{([a]+1)M-1}|F(T^ky)|=[a]Av(|F|)$ which is
bounded due to the $S$-integrability of $F$.

Let $B=\{T^ky\ |\ k=L,\dots,K-1\}$. Then
$\frac{|B|}M=\frac{K-L}M\approx 0$. Thus, $\frac
1M\sum\limits_{y\in B}|F(y)|\approx 0$, due to the
$S$-integrability of $F$.
So, $V=\frac MK\cdot\frac 1M\sum\limits_{y\in B}|F(y)|\approx 0$ $\Box$

The following example shows that this quasi-proposition may fail for the case of very big $K,L$ such that
$\frac KM\approx\frac LM\approx 0$.

\textbf{Example 3} Let $Y$ and $T$  be the same as in the previous examples. Fix a number $K\sim\infty$ such that $K/M\approx 0$ and consider the function $F:Y\to\*\R$ given by the formula:

\begin{equation}\label{ex03}
F(k)=\left\{\begin{array}{ll} 1,\ &mK\leq k<(m+1)K,\ m<R,\ m\ \mbox{is even} \\
0,\ &(mK\leq k<(m+1)K,\ m<R,\ m\ \mbox{is odd})\lor RK\leq m<M
\end{array}\right.
\end{equation}

The function $F$ is bounded and, thus, $S$-integrable.

To show that Theorem \ref{ErgMeanStab} fails for the set
$\{L\sim\infty\ |\ \frac LM\approx 0\}$ it is enough to prove that
$A_K(F,T,y)\not\approx A_{\left[\frac K2\right]}(F,T,y)$ for all $y$
in some set of the positive measure $\mu$

Let $D=\bigcup\limits_{m<R}\{y\in Y\ | mK\leq k<mK+\frac K2\}$. Then, $\mu_L(D)=\frac
12$. It is easy to see, that $\all\, y\in D\, \all n\leq\frac K2\
F(T^n(y))=F(y)$, thus, $A_{\left[\frac K2\right]}(F,T,y)=G(y)$.

For every number $n\ll\infty$ consider the set
$D_n=\{y\in D\ |\ |A_K(F,T,y)-F(y)|<\frac 1n\}$. It is enough to
prove that $\lim\limits_{n\to\infty}\mu_L(D_n)=0$. Since the
cardinality of the set $D_n\cap [mK,(m+1)K)$ is the same for all
$m<R$, it is enough to calculate the cardinality of $E_n=D_n\cap
\left[0,\frac K2\right]$.  Since $E_n\ss D$, for any $y\in E_n$ one has
$A_{\left[\frac K2\right]}(F,T,y)=F(y)=1$. On the other hand
$A_K(F,T,y)=\frac{K-y}K=1-\frac yK$. So, $|E_n|=\left[\frac
Kn\right]$ and $\mu_L(D_n)=\left(R\left[\frac
Kn\right]\right)/(RK+S)\to 0$. $\Box$

\begin{figure}
\centering
\includegraphics*[viewport= 50 500 500 750, scale=0.5]{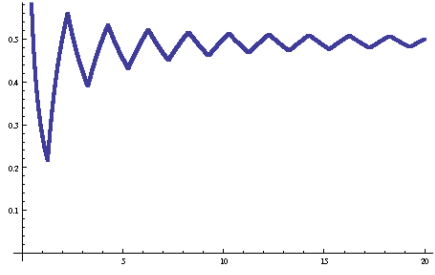}
\caption{$\Gamma(F)$ for $F$ defined in Eq.\ref{ex03}}
\label{fig3}
\end{figure}

On Fig.\ref{fig3} we see the visual image of the set $\Gamma(F)$ for $M=100000$, $K=1000$ and the randomly chosen $y=870722$.
\begin{Th} \label{NSBET} If $Y$ is a hyperfinite set, $|Y|=M$ is a very big number, and $T:Y\to Y$ is an arbitrary permutation, then for any $S$-integrable function $F:Y\to\*\R$ there exists an initial segment $I\ss\iN$  such that for almost all $y\in Y$ for all $L,N\in I$ one has $A_L(F,T,y)\approx\A_N(F,T,y)$. In the case of transitive permutation $T$ one has
\begin{equation} \label{NSBET1}
A_N(F,T,y)\approx \o Av(F)=\int_Y\o Fd\mu_L
\end{equation}
for all $N$ such that $\frac NM\approx 1$.
\end{Th}

The statement of this theorem cannot be seen on the on the picture of $\Gamma(F)$, where we see simply some ambiguity around the origin. Theorem \ref{NSBET}  states, however, that one can always find  a number $N\sim\infty$ such that the visual image of the set of points $\{(n/N,A_n(F,T,y))\ |\ n=1,...,N\}$, is a horizontal line (maybe again with some ambiguity around the origin due to numbers $n\ll\infty$) for almost all $y\in Y$. In other words, this means that it is always possible to find a microscope, such that looking through it on the picture of $\Gamma(F)$ around the origin one can always see that the initial part of this graph is a horizontal line
and we may use one microscope for almost all $y\in Y$

\begin{figure}
\includegraphics*[viewport=50 500 500 750,scale=0.4]{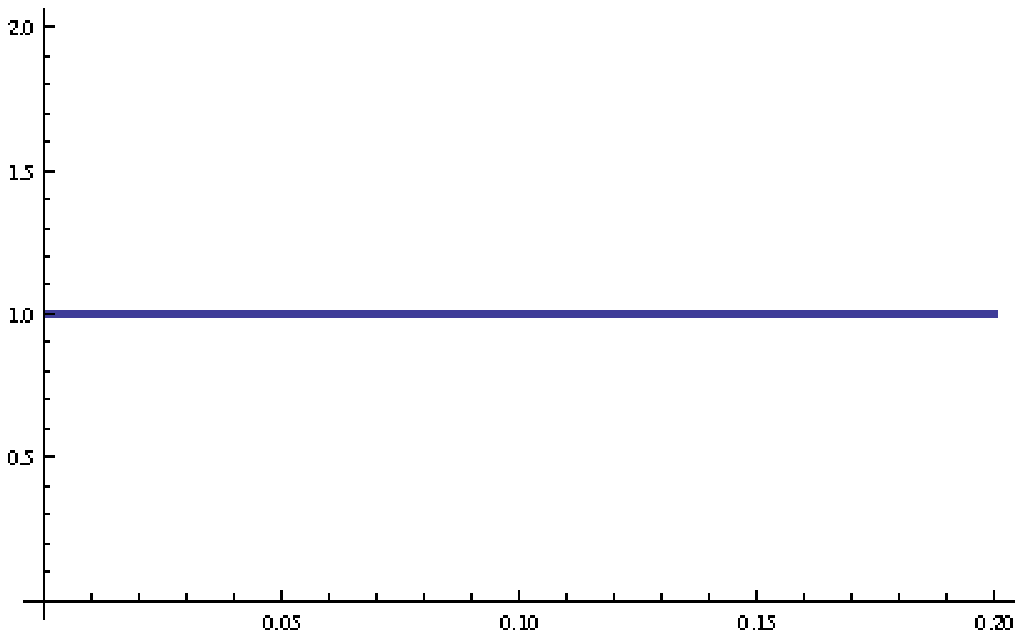}
\hfill
\includegraphics*[viewport=50 500 500 750,scale=0.4]{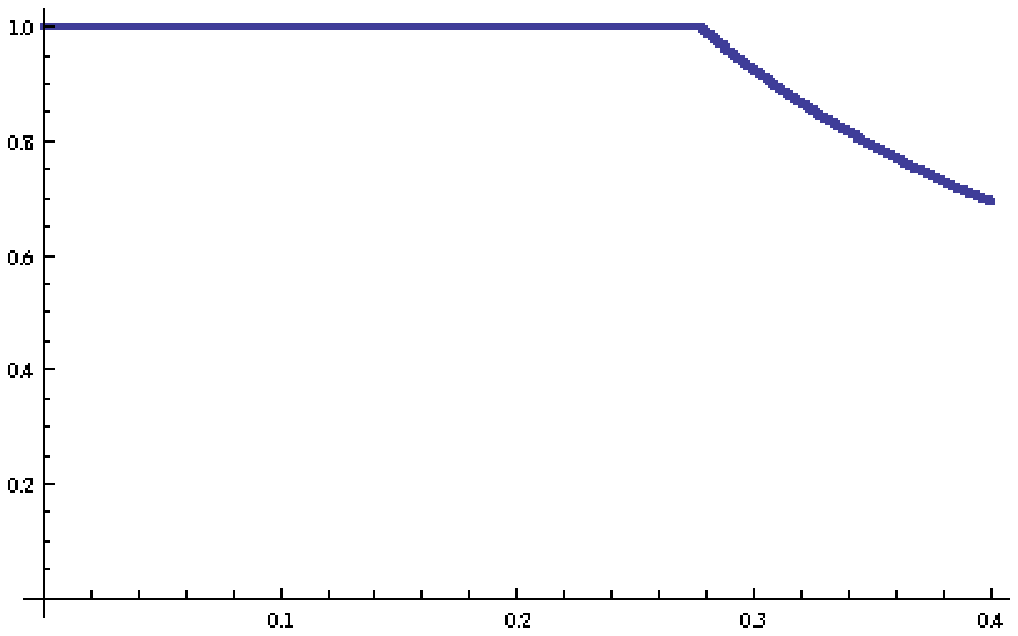}
\caption{$\Gamma(F)$}
\label{fig4}
\end{figure}

On  Fig.\ref{fig4} one can see the visual image of the set $\{(n/N,A_n(F,T,y))\ |\ n=1,...,N\}$ for the function $F$
of Example 3 where $M=100000$, $K=1000$, $N=0.2K$ for. Fig.4 shows, that
$N=0.4K$ does not satisfy Theorem \ref{NSBET}.

The proof of Theorem \ref{NSBET} is not elementary. It is contained in Section 4. It is interesting that this proof uses Ergodic Theorem and Egoroff's Theorem for infinite probability spaces. One can easily show that any $N\ll K$ in Example 3 satisfies Theorem \ref{NSBET} .

The following example shows that the statement Theorem of \ref{NSBET} may not be true for all $y\in Y$.

\textbf{Example 4.} Let $Y,T$ be the same as in the previous examples and $F(y)=\frac
yM$. Consider the following function $\psi:[0,1]^2\to\R$
by the formula
$$
\psi(a,t)=\left\{\begin{array}{l} t+\frac a2,\ 0\leq t\leq 1-a
\\ t+\frac a2-1+\frac 1a(1-t),\ 1-a<t\leq 1\end{array}\right.
$$
Then it is easy to calculate, that
\begin{equation} \label{examle-4}
A_K(F,T,y)\approx\left\{\begin{array}{l} \psi(0,\frac yM),\
\frac KM\approx 0,\ \forall y<M-K \\
\ \psi(a,\frac yM)\ \frac KM\approx a>0,\
\forall y\in Y
\end{array}\right.
\end{equation}

Another easy calculation shows that for $K\sim\infty$ such that $\frac KM\approx 0$, and for $y\in\{K-M,...M-1\}$ one has $A_K(F,T,y)\approx\infty$. Obviously $\mu(\{K-M,...M-1\})=\frac{M-K}M\approx 0$.

\bigskip

iii) In this section we introduce definition of approximation of a dynamical system $(X,\nu\tau)$ on a compact metric space $(X,\rho)$ with a Borel measure $\nu$ and
a $\nu$-preserving transformation $\tau$ by an internal dynamical system $(Y,\mu,T)$ on a hyperfinite set $Y$, $|Y|=M$ with a uniform probability measure $\mu$ and an internal permutation $T:X\to X$. Since we consider only the measure $\mu$ on a hyperfinite space we write $(Y,T)$ for the above hyperfinite dynamical system. For a set
$C\ss X$ we denote the set $\{x\in X\ |\ \ex\, c\in C\ \rho(x,c)<\e\}$ by $U_{\e}(C)$.

\begin{Def} \label{hypap} 1) Let $\f:Y\to \* X$ be an internal injective map such that
for every closed set $C\ss X$ there exists an initial segment $I\ss\*\N$ such that $\mu(\f^{-1}(U_{\frac 1N}(\*C)))\approx\nu(C)$ for all $n\in I$.
Then the pair $(Y,\f)$ is said to be a \emph{hyperfinite approximation (h.a.)} of the measure space $(X,\nu)$.
In case of $Y\ss\*X$ and the identical embedding $\f$ we say that $Y$ is a h.a. of $(X,\nu)$. Obviously, any h.a. $(Y,\f)$ is equivalent to the h.a. $\f(Y)$.

2) Let $\tau:X\to X$ be a measure preserving transformation of $X$ and $(Y,\f)$ be a h.a. of $(X,\nu)$. Then we say that an internal permutation
$T:X\to X$ is a h.a. of the transformation $\tau$ if for almost all $y\in Y$ one has $\f(T(y))\approx\tau(\f(y))$.
We say also that the internal triple $(Y,T,\phi)$ is a h.a. of the dynamical system 
$(X,\nu,\tau)$.
\end{Def}

\begin{Prop} \label{pres-mes} A pair $(Y,\f)$ is h.a. approximation of $(X,\nu)$, if the map $\sp\circ\f\:Y\to X$ is a measure preserving map with respect to the measure $\nu$ and the Loeb measure $\mu_L$. \end{Prop}

\textbf{Proof} By Proposition \ref{st-nonst} the condition of Definition \ref{hypap} (1) is equivalent to the following condition
\begin{equation} \label{hypap1}
\nu(C)=\lim_{\e\to 0}\o\mu(\f^{-1}(\*U_{\e}(C))).
\end{equation}
It is easy to see that, if $C$ is a compact set, then 
$\bigcap\limits_{n\in\iN}U_{\frac 1n}(\*C)=\sp^{-1}(C)$. Using this fact and the equality (\ref{hypap1}) one obtains the equality $\nu(C)=\mu_L(\f^{-1}(\sp^{-1}(C)))$. $\Box$

Let $f\in L_1(\nu)$ and $Y\ss\*X$ be a h.a. of $(X,\nu)$. Restrict the function $\sp:\*X\to X$ on the set $Y$. Then Proposition \ref{pres-mes} shows that $\sp$ is a measure-preserving map and, thus, $f\circ\sp\in L_1(\mu_L)$. Due to Proposition \ref{S-int}(3) $\f\circ\sp$ has an $S$-integrable lifting $F$. We say in this case that $F$ is an $S$-integrable lifting of the function $f$. In this case
\begin{equation} \label{hypap2}
Av(F)=\int_YFd\mu_L=\int_Xfd\nu.
\end{equation}

\begin{Prop} \label{cont-lif}
If $f\in C(X)$, then $\*f\upharpoonright Y$ is an $S$-integrable lifting of $f$.
\end{Prop}

\begin{Th} \label{cycle-appr} Let $\nu$ be a non-atomic Borel measure on a compact metric space $X$ such that the measure of every ball
is positive. Then
\begin{enumerate}
\item for every set $A\ss X$ such that $\nu(A)=1$ there exists a hyperfinite set $Y\ss \*A$ such that $Y$ is a
h.a. of $(X,\nu)$.
\item For every dynamical system $(X,\nu,\tau)$  and for every h.a. $Y$ of $(X,\nu)$ there exists a h.a. $(Y,T)$.
Moreover, one can choose a h.a. $T$ of $\tau$ to be a cycle of the length $M$.
\end{enumerate}
\end{Th}
The proof of this theorem is contained in Section 4.

In what follows a $Y$-cycle of the length $M$ is said to be a \emph{transitive permutation} of $Y$.

\textbf{Example 4'.} Returning to Example 4 above define the map $\f:Y\to [0,1]$ by the formula $\f(y)=\frac yM$. Then, obviously, the pair $(Y,\f)$ approximates the probability space $([0,1],dx)$, where $dx$ is the Lebesgue measure. The permutation $T$ approximates the identical map $id:[0,1]\to [0,1]$. Indeed, for all $y\neq M-1$ one has $id(\f(y))=j(y)=\frac yM\approx \frac{y+1}M=\f(T(y))$. Thus, $Y,\ \f$, $T$ and $\tau=id$ satisfy Definition \ref{hypap}. The function $F$ of Example 4 is a lifting of the function $g(x)=x$

 \begin{Prop} \label{CNSBET} Let $(Y,T)$ be a h.a. of $(X,\nu,\tau)$, $f\in L_1(\nu)$,  $\tilde{f}=\lim\limits_{n\to\infty}A_n(f,\tau, x)$ and
let $\widetilde{F}$ be an $S$-integrable lifting of $\tilde{f}$, then there exists an initial segment $I\ss\iN$,
such that for  almost all $y\in Y$
$$\all\, K\in I\ A_N(F,T,y)\approx \widetilde{F}(y)\approx \widetilde f(\o y).$$.
\end{Prop}

This proposition follows immediately from Theorem \ref{NSBET}.

\begin{Cor} \label{non-erg} Let $T$ be a transitive permutation and let $\tau$
be a non-ergodic transformation. Consider a function
$f\in L_1(\nu)$ such that the set $B\ss X$ of all $x\in X$ satisfying
inequality $\lim\limits_{n\to\infty}A_n(f,\tau,x)\neq Av(f)$ has a
positive measure $\nu$. Then there exist infinite
$M$-bounded $N,K$ such that for almost all $y\in\sp^{-1}(B)$ one
has $A_N(F,T,y)\not\approx A_K(F,T,y)$.
\end{Cor}

\textbf{Proof} \newcommand{\wt}{\widetilde}. Let $\wt
f=\lim\limits_{n\to\infty}A_n(f, \tau,\cdot)$ and $\wt F$ be the same
as in Proposition \ref{CNSBET}. By this proposition there exists $N\in\iN$
such that $\frac{N}{M}\approx 0$ and $A_N(F,T,y)\approx\widetilde{F}(y)$ $\mu_L$-a.e.
Thus, $A_N(F,T,y)\not\approx Av(f)$ for $\mu_L$-almost all $y\in\sp^{-1}(B)$.

On the other hand, since $T$ is a cycle of length $M$, by Theorem \ref{ErgMeanStab} one has
$A_K(F,T,y)\approx Av(F)\approx Av(f)$ for all $y\in Y$ and for all $K$ such that $\frac{K}{M}\approx 1$. Thus,
$A_K(F,T,y)\not\approx A_M(F,T,y)$ for $\mu_L$-almost all $y\in\sp^{-1}(B)$.
$\Box$

Proposition shows that if $(Y,T)$ is a h.a. of a dynamical system $(X,\nu,\tau)$, $\tau$ is a non-ergodic transformation and $T$ is a transitive permutation, then
there exists a function $f\in L_1(\nu)$ an internal set $B\ss Y$, $\mu(B)\gg 0$ and $K,L\in\iN$ such that $A_K(F,T,y)\not\approx A_L(F,T,y)$ for all $y\in B$, where
$F$ is an $S$-integrable lifting of $f$.

\begin{Cor}\label{non-erg2} If under the conditions of the previous paragraph, for any $f\in L_1(\mu)$ for almost all $y\in Y$ and for all $N\in\iN$ one has
\begin{equation}\label{appr-eq}
A_N(F,T,y)\approx\int_Xfd\nu,
\end{equation}
then $\tau$ is an ergodic transformation\end{Cor}.

We do not know, wether the sufficient condition of the ergodicity of $\tau$ formulated in Corollary \ref{non-erg2}, is also a necessary condition. By Proposition \ref{CNSBET} the approximate equality (\ref{appr-eq}) holds for all $N$ in some initial segment $I\ss\iN$ for almost all $y\in Y$ for an ergodic transformation $\tau$ or
for any transformation $\tau$ and transitive permutation $T$, if $\frac NM\approx n\in\N$ or, if $\frac NM\sim\infty$ (see Theorem \ref{INTBET}).

This and even stronger is necessary for uniquely ergodic transformations. Recall that a continuous transformation $\tau:X\to X$ is said to
be \emph{uniquely ergodic} if there exists only one $\tau$-invariant Borel measure on $X$\footnote{Krylov-Bogoljubov theorem claims the existence of at least one $\tau$-invariant measure.}.

\begin{Th} \label{un-erg} If $\tau$ is a uniquely ergodic transformation of a compact metric space $X$,
$Y\ss \*X$ is a hyperfinite dense subset of $\*X$, and $T:Y\to Y$ is an internal permutation such that $\all\,y\in Y\ \sp (T(y))=\tau(\sp (y))$,
then for every $y\in Y$ such that the $\tau$-orbit of $\sp(y)$ is dense in $X$, for every $N\in\iN$ and for every $f\in C(X)$ one has
\begin{equation} \label{un-erg-1}
A_N(\*f\upharpoonright Y,T,y)\approx\int\limits_Xfd\nu,
\end{equation}
where $\nu$ is the $\tau$-invariant measure.
\end{Th}

The proof of this theorem is contained in Section 4.

\textbf{Example 5.} In this example we consider a hyperfinite set $Y=\{0,1,...,M-1\}$ and a permutation $T:Y\to Y$ given by the formula $T(y)=y+P\mod M$. We choose $P$ and $M$ to be relatively prime, so that $T$ is a cycle of length $M$. The approximation $(Y,\f)$ of the probability space $([0,1], dx)$ used in Example 4' can be used as well for the probability space $([0,1), dx)$, where $[0,1)$ is provided by the topology of the circle. For any $t\in\R$ the measure preserving transformation $\tau:[0,1)\to [0,1)$, such that $\tau(y)=y+t(\mod 1)$ is called the $t$-shift of a circle. This transformation is continuous on $[0,1)$. It is easy to see that if $\frac PM\approx t$, then the triple $(Y,\f,T)$ approximate the shift
$\tau$. We present the visual images of $\Gamma(F)$ for $F=f\upharpoonright \f(Y)$, where
$$f(x)=\left\{\begin{array}{ll}\frac{10x}9,\ &\mbox{if}\ 0\leq x<0.9\\ 10(1-x),\ &\mbox{if}\ 0.9\leq x<1\end{array}\right.$$

\begin{figure}
\includegraphics*[viewport=50 500 500 750, scale=0.4]{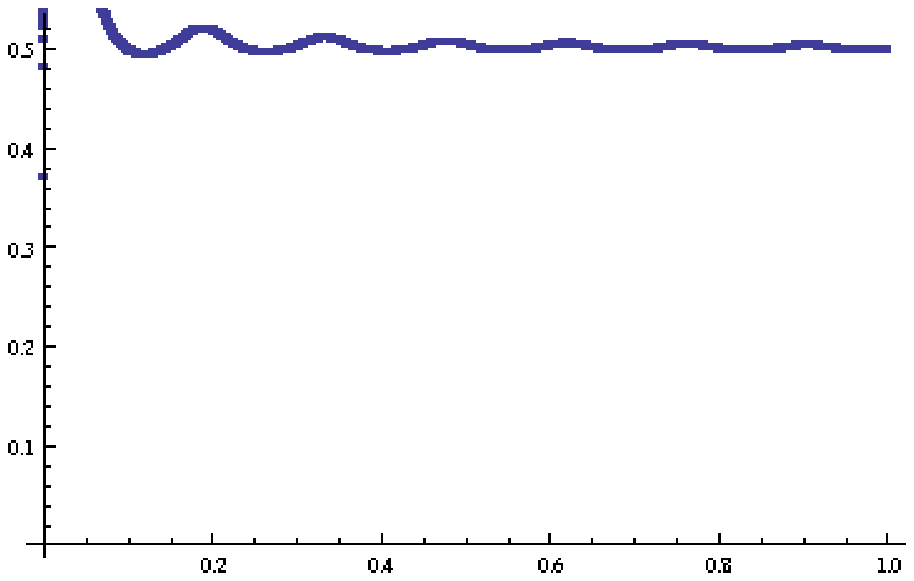}
\hfill
\includegraphics*[viewport=50 500 500 750, scale=0.4]{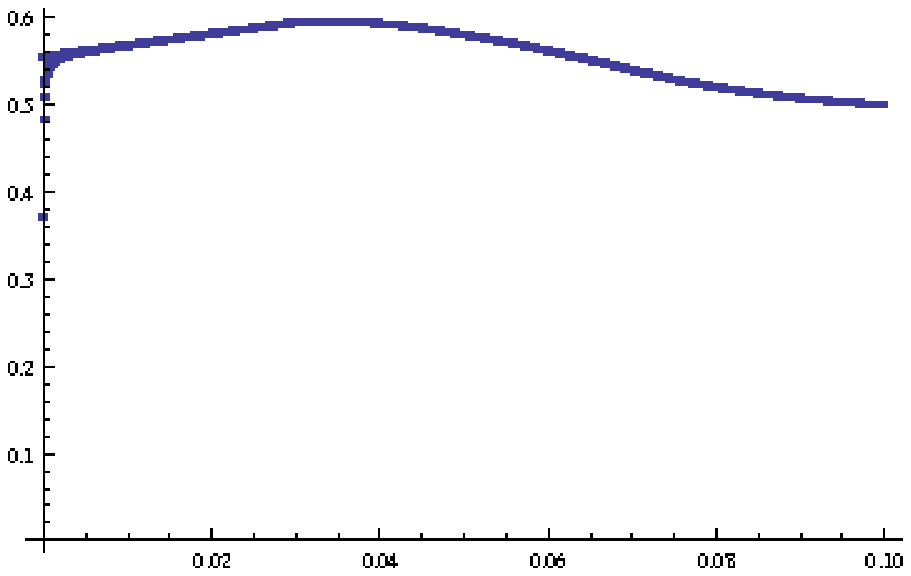}
\caption{$\Gamma(F)$ for transformation $\tau(x)=x+2/3 \mod 1$}
\label{fig5}
\end{figure}

We choose the function $f$ that is close to the function $g(x)=x$, considered in Example 4'. However, $f$ is continuous on the circle $[0,1)$, while $g$ is discontinuous at $x=0$. We consider two cases.

a). On Fig.\ref{fig5} the visual image of $\Gamma(F)$ on the interval $[0,1]$ and at the neighborhood of $0$ are presented for the case of $M=33334$, $P=22225$ and the randomly chosen $y=16667$.

\begin{figure}
\includegraphics[viewport=50 500 500 800, scale=0.4]{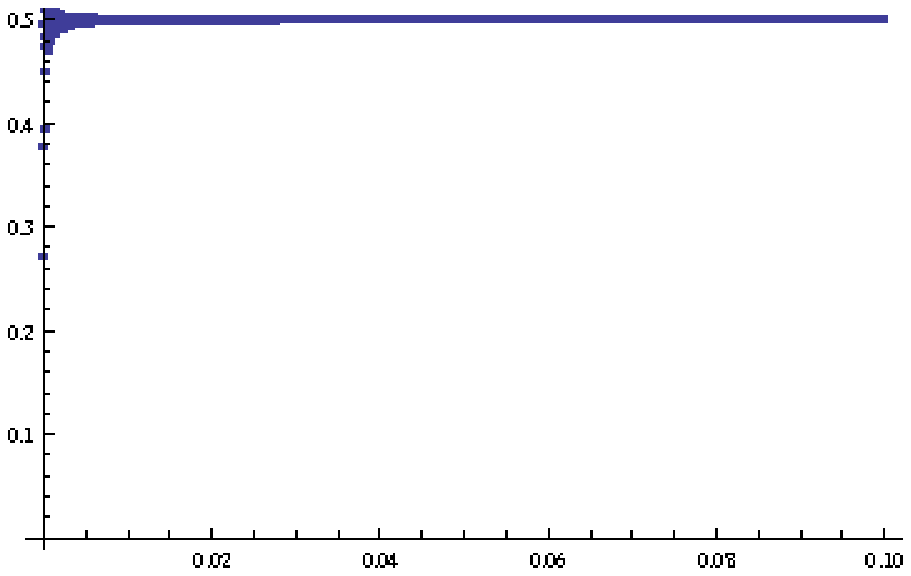}
\hfill
\includegraphics[viewport=50 500 500 800, scale=0.4]{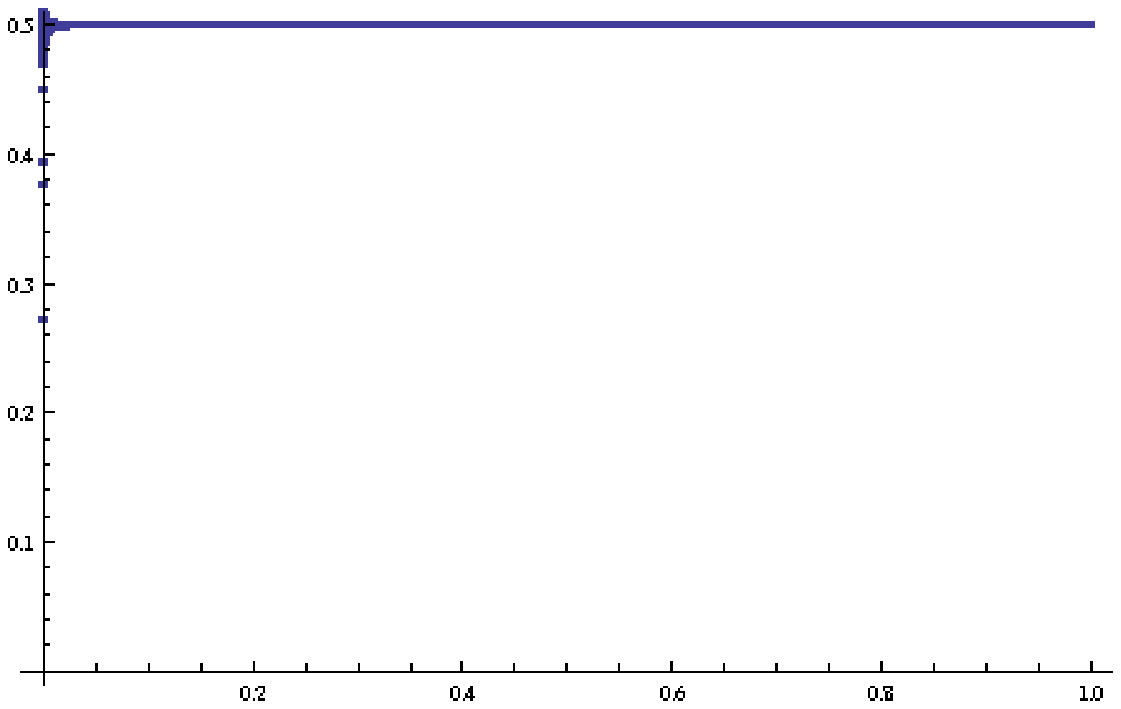}
\caption{$\Gamma(f)$ for transformation $\tau(x)=x+1/\sqrt{2} \mod 1$}
\label{fig6}
\end{figure}

b). On Fig.\ref{fig6} we see the visual image of $\Gamma(f)$ on the interval $[0,1]$ and at the neighborhood of $0$ are presented for the case of $M=25001$, $P=17677$ and the randomly chosen $y=6119$. In this case $(Y,\f,T)$ approximates $1/\sqrt{2}$-shift.
In accordance with Theorem \ref{un-erg} this visual image is a horizontal line.

To explain the difference in the graphs in the cases a) and b) notice that in the case a) $\frac PM\approx \frac 23$. Actually, $\left|\frac PM-\frac 23\right|\leq 0.00046$ that is enough for our problem to consider these numbers to be infinitesimally close (see the discussion in the Example 1). In this case $x=\frac yM=\frac{16667}{33334}\approx0.5$ and $\tau(y)=y+\frac 23 (\mod 1)$. It is easy to see that for $\frac 13<y<\frac 23$ one has $\hat{\f}(y)=\frac 13[\f(y-1/3)+\f(y)+\f(y+1/3)]$ (for any integrable function $\f$). So, in our case $\hat{\f}(0.5)\approx 0.56$.
We see that the value of the function on Fig.5 at the neighborhood of $0$ is close to $0.56$ that agrees with Proposition \ref{CNSBET}.
Since $T$ is a cycle of length $M$ the value of the function on Fig.5 at the neighborhood of $1$ is close to $Av(f)=Av(\f\upharpoonright j(X))\approx\int_0^1\f d\mu$ by Definition \ref{hypap}. In our case $\int_0^1\f dx=0.5$. The visual image of $\Gamma(f)$ is a continuous function on $[0,1]$ in accordance with Theorem \ref{ErgMeanStab}.

In case b) $\frac PM\approx\frac 1{\sqrt 2}$ ($|\frac PM-\frac 1{\sqrt 2}|\leq 0.00006$). So, in this case $T$ approximates irrational shift of the circle $[0,1)$. It is well-known that irrational shifts of the circle are uniquely ergodic transformations. Since 
$\f$ is a continuous function, $A_K(f,T,x)\approx 0.5$ for all $K\sim\infty$ by Theorem \ref{un-erg}. Thus, the visual image of $\Gamma(f)$ is the horizontal line $y=0.5$, that is perfectly reflected on both pictures of Fig.\ref{fig6}.

The consideration of these two examples arises the following question. Suppose that we have a ratio $\frac PM$ of two relatively prime numbers. In what case this ratio can be considered as "practically" rational number and in what case one should consider it as "practically" irrational number.  Using the informal language of the Introduction one can say that that $\frac PM$ is "practically" rational, if there exist two natural numbers $m,n\ll\infty$ such that $\frac PQ\approx \frac mn$, and $\frac PM$ is "practically" irrational otherwise. Certainly the exact answer strongly depends on a problem, in which we need to answer this question. However, the statements of this section provide us with some qualitative understanding of the correlation between the behaviors of very big discrete systems and their continuous approximations, (or vise versa, continuous systems and their discrete approximations).

{\bf Example 6}. (Approximations of Bernoulli shifts). Let $\Sigma_m=\{0,1,\dots m-1\}$. Consider the compact space $X=\Sigma_m^{\mathbb{Z}}$ with the Tychonoff topology. Let $a$ be a function, such that $\dom(a)\subset\Z$ is finite, and $\range(a)\ss\Sigma_m$. Let $S_a=\{f\in X\ |\ f\upharpoonright\dom(a)=a\}$. Then the family of all such $S_a$ form a base of neighborhoods of the compact space $X$. For $g\in\*X$ set $f=g\upharpoonright\Z$, then $f\in X$ and it is easy to see that $f=\sp(g)$.

The continuous transformation $\tau:X\to X$ defined by the formula $\tau(f)(n)=f(n+1)$ where $f\in X$ and $n\in\Z$ is an invertible Bernoulli shift.
Every probability distribution $\{p_0,\dots,p_{m-1}\}$ ($p_i>0,\ \sum\limits_{i=0}^{m-1}p_i=1$) on $\Sigma_m$ defines a Borel measure on $X$ that is obviously
invariant with respect to $\tau$. It is well-known that $\tau$ is ergodic for each of these measures. So, the transformation $\tau$ is not uniquely ergodic.
Here we restrict ourselves only to the case of the uniform distribution on $\Sigma_m$, i.e. to the case of $p_0=\dots=p_{m-1}=\frac 1m$.
The corresponding Borel measure on $X$ is denoted by $\nu$. Obviously, $\nu(S_a)=m^{-|\dom(a)|}$.

We construct here two hyperfinite approximations of the dynamical system $(X,\nu,\tau)$. First we consider the straightforward approximation by a hyperfinite shift.
Fix $N\in\iN$ and set $Y=\Sigma_m^{\{-N,\dots,N\}}$. Then $M=|Y|=m^{2N+1}$. Define $\l:Y\to\*X$ as follows. For $y\in\Sigma_m^{\{-N,\dots,N\}}$ set
\begin{equation} \label{3}
\l(y)(n)=\left\{\begin{array}{ll} y(n),\ &|n|\leq N \\ 0,\ &|n|>N\end{array}\right.
\end{equation}
Then $\sp\circ\l (y)=y\upharpoonright\Z$. Thus, for every standard neighborhood $S_a$
defined above one has $(st\circ\l)^{-1}(\*S_a)=\{y\in Y\ |\ y\upharpoonright\dom(a)=a\}$.
So, $\mu_L((\sp\circ\l)^{-1}(\*S_a))=\nu(S_a)=m^{-\dom(a)}$.
This proves that $(Y,\l)$ is a h.a. of $(X,\nu)$.

Certainly, an arbitrary internal map from $Y$ to $\Sigma_m^{\*\Z\setminus\{-N,\dots,N\}}$
can be used to define the values $\l(y)(n)$ for $|n|>N$ and $y\in Y$ in the definition of $\l$ (\ref{3})

In what follows we use notations $y_1\approx y_2$ and $\sp(y)$ for $\l(y_1)\approx\l(y_2)$ and $\sp(\l(y))$ respectively.

Define the map $S\to S$ by the formula $S(y)(n)=y\left(n+1(\mod 2N+1)\right)$ for any $y\in Y$ and $n\in\{-N,\dots,N\}$. Then obviously $\tau(\sp(y))=\sp(S(y))$ for all $y\in Y$. So, $(Y,\l,S)$ is a h.a. of the dynamical system $(X,\nu,\tau)$.

Since every point $y\in Y$ is $(2N+1)$-periodic with respect to $S$ the permutation $S$ is not transitive. Though the existence of a transitive h.a. of $\tau$ is proved in Theorem \ref{cycle-appr}, it is not easy to construct such approximation explicitly.

To do this we reproduce here the construction of de Bruijn sequences.

\begin{Def} \label{de Bruijn}
An $(m,n)$-de Bruijn sequence on the alphabet $\Sigma_m$ is a sequence $s=(s_0,s_1,\dots,s_{L-1})$ of $L=m^n$ elements $s_i\in\Sigma_m$ such
that all consecutive subsequences $(s_i, s_{i\oplus 1},\dots,s_{i\oplus n-1})$ of length $n$ are distinct.

Here and below in this example the symbols $\oplus$  and $\ominus$ denote + and -  modulo $L$, so that the sequence $s$ is considered as a sequence of symbols from
$\Sigma_n$ placed on a circle.
\end{Def}

It was proved \cite{deB1, deB2} that there exist $(m!)^{m^{n-1}}\cdot m^{-n}$ $(m,n)$-de Bruijn sequences. See also \cite{deB3} for a simple algorithm for de Bruijn sequences and more recent references.

To construct a transitive h.a. $T:Y\to Y$ of $\tau$ fix arbitrary $(m, 2N+1)$ de Bruijn sequence $s=(s_0,s_1,\dots,s_{M-1})$ here $L=M$.
Let $y=(y_{-N},\dots,y_{-N})\in Y$. Then there exists the unique consecutive subsequence $\s(y)=(s_i,s_{i\oplus 1},\dots,s_{i\oplus 2N})$ such that
$y_j=s_{j\oplus i\oplus N})$.  Set $P(\s(y))=(s_{i\oplus 1},\dots, s_{i\oplus 2N\oplus 1})$ and $T(y)=\s^{-1}(P(\s(y)))$. Notice that if $i< M-N$ , then
for all $j\leq N$ one has $s_{j\oplus i\oplus N}=s_{j+i+N}$. So, $T(y)_j=y_{j+1}$ for all $j<N$ and, thus, for all standard $j$.
This last equality implies that $\sp(T(y))=\tau(\sp(y))$ for all $y\in Y$ such that the first entry of the sequence $\s(y)$ is the $i$-s term of the initial de Bruijn
sequence for $i \leq L-2N-1$. So, $\mu_L(\{y\ |\ \sp(T(y))=\tau(\sp(y))\})\geq\frac{M-N}M\approx 1$. This proves that $T$ is a h.a. of $\tau$. We call $T$ a de Bruijn approximation of $\tau$.

It is interesting to study the behavior of ergodic means of described approximations. This problem will be discussed in another paper. We confine ourselves with two simple remarks.

1. If $\s(y)=\la s_i,\dots,s_{i+2N}\ra$ and $i<M-N$, then $A_n(F,T,y)=A_n(F,S,y)$ for all 
$n<N$.

2. Let $S_0=\{f\in X\ |\ f(0)=1\}$, so that $\nu(S_0)=\frac 12$ and let $\chi_0$ be 
a characteristic function of $S_0$. For $y\in Y$ let $f=\sp(y)$. Set 
$A(y)=f^{-1}(\{1\})\cap\N$. Recall that the density of $A(y)$ is given by the formula
$$
d(A(y))=\lim\limits_{m\to\infty}\frac{|A(y)\cap\{0,\dots, m-1\}|}m.
$$
It is easy to see that for $m<N$ one has
$$A_m(\*\chi_0,T,y)=\frac{|A(y)\cap\{0,\dots, m-1\}|}m.$$
So, for all $y\in Y$ such that the density $d(A(y))$ exists one has 
$\ex\, K\in\iN\,\all\, m\in\iN\;\;\;\; m\leq K\Limpl A_m(\chi_0,T,y)\approx d(A(y))$.
Due to Proposition \ref{CNSBET} there exist $K\in\iN$ such that for $\mu_L$-almost all $y\in Y$ one has $A_m(\*\chi_0, T,y)\approx\frac 12$.

\bigskip

iv) (\emph{Formulation of results in the framework of the standard mathematics.}) While in nonstandard analysis we use the notion of an infinite number (hyperfinite set) as a formalization of the notion of a very big number (finite set), in classical mathematics we use the sequences of numbers (finite sets) diverging to infinity to formalize these notions. For example, in previous sections we considered a hyperfinite set $Y$ and its internal permutation $T:Y\to Y$. If we want to treat the same problems in the framework of standard mathematics, we have to consider a sequence $(Y_n, T_n)$ of finite sets $Y_n$ whose cardinalities tend to infinity and their permutations $T_n$. Similarly, internal functions $F:Y\to\*\R$ correspond to sequences $F_n:Y_n\to\R$ in standard mathematics.

First, we discuss what property of such sequences correspond to the property of an internal function $F$ to be $S$-integrable. The following proposition gives a
reasonable answer to this question.

\begin{Prop} \label{uni-int} Let $Y_n$ be a standard sequence of finite sets, such that $|Y_n|=M_n\to\infty$ as $n\to\infty$. Then for an arbitrary sequence
$F_n:X_n\to\R$ the following statements are equivalent:
\begin{enumerate}
\item For every $K\in\iN$ the function $\* F_K$ is $S$-integrable.
\item
\begin{equation}\label{uni-int-1}
\lim\limits_{n,k\to\infty}\frac 1{M_n}\sum\limits_{\{x\in Y_n |\ |F_n(x)|>k\}}|F_n(x)|=0
\end{equation}
\end{enumerate}
\end{Prop}

The proof can be obtained easily from the definition of $S$-integrable functions (\ref{S-int-1}) using arguments similar to those that were used in the proof of
Theorem \ref{sequences}

A sequence $F_n$ that satisfies the statement (2) of Proposition \ref{uni-int} is said to be \emph{uniformly integrable}.

Proposition \ref{uni-int} leads to establishing the standard version of Theorem \ref{ErgMeanStab}.

\begin{Prop} \label{ErgMeanStabSt} In conditions of Proposition \ref{uni-int} let $T_n:Y_n\to Y_n$ be a sequence of transitive permutations and $F_n:Y_n\to\R$ be a uniformly integrable sequence. Consider two sequences of natural numbers $K_n$ and $L_n$ such that $\frac{K_n}{M_n}$ is bounded, $\liminf\frac{K_n}{M_n}>0$ and $\lim\limits_{n\to\infty}\frac{K_n}{L_n}=1$. Then the following two statements  are true.
\begin{enumerate}
\item For any $\e>0$ one has
\begin{equation}\label{ErgMeanStab 0}
\lim\limits_{n\to\infty}\frac 1{M_n}\cdot\left|\{y\in Y_n\ |\ \left|A_{K_n}(F_n,T_n,y)-A_{L_n}(F_n,T_n,y)\right|\geq\e\}\right|=0
\end{equation}
\item. If $T_n$ is a sequence of transitive permutations or $F_n$ is a sequence of bounded functions, then
\begin{equation}\label{ErgMeanStab1}
\lim\limits_{n\to\infty}\max\limits_{y\in Y_n}\left|A_{K_n}(F_n,T_n,y)-A_{L_n}(F_n,T_n,y)\right|=0
\end{equation}
\end{enumerate}
\end{Prop}

 It is not too difficult to rewrite the proof of Theorem \ref{ErgMeanStab} in (standard) terms of Proposition \ref{ErgMeanStabSt}. It is also easy to deduce Proposition \ref{ErgMeanStabSt} from Theorem \ref{ErgMeanStab} using arguments close to those of the proof of Theorem \ref{sequences}.

The rigorous interpretation of Theorem \ref{NSBET} in the framework of standard mathematics is much more difficult.  To formulate the corresponding rigorous theorem, we need to remind the notion of an ultrafilter and of the limit of a sequence over a non-principle ultrafilter.

Recall that a subset $\F\ss\P(\N)$ is said to be a \emph{non-principle ultrafilter}, if $\F=\{A\ss\N\ |\ m(A)=1\}$ for some finitely additive measure $m$ on $\P(\N)$ that assumes only two values $0$ and $1$  such that $m(B)=0$ for any finite set $B\ss\N$ and $m(\N)=1$. For a sequence $\a:\N\to\R$ we say that $\lim\limits_{\F}\a=L$, if for any $\e>0$ the set $\{n\in\N\ |\ |\a(n)-L|<\e\}\in\F$. It is known that any bounded sequence has a limit over any non-principle ultrafilter. For two sequences $\a,\b:\N\to\R$ we say that $\a\leq_{\F}\b$ iff $\{n\in\N\ |\ \a(n)\leq\b(n)\}\in\F$.

\begin{Prop}\label{NSBETst} Let $T_n:Y_n\to Y_n$ be a sequence of arbitrary permutations. Then for every non-principle
ultrafilter $\F\ss\P(\N)$ there exists a sequence $K_n\to\infty$ as $n\to\infty$ and a sequence $A_n\ss Y_n$ such that $\lim\limits_{\F}\frac{|A_n|}{M_n}=0$ and for any $L_n\to\infty$ as $n\to\infty$, if $\la L_n\ra\leq_{\F}\la K_n\ra$, then $\lim\limits_{\F}(A_{K_n}(f_n,T_n,x_n)-A_{L_n}(f_n,T_n,x_n))=0$.\end{Prop}

This proposition doesn't have such intuitively clear sense as Theorem \ref{NSBET}. One hardly can find a proof of it, that doesn't use the ideas of nonstandard analysis.

To formulate the standard version of Definition \ref{hypap} introduce the following notation.
Let $Z\ss X$ be a finite subset of $X$ and $\d_Z=\frac 1{|Z|}\sum\limits_{z\in Z}\d_z$, where
$\d_Z$ is a Dirac measure at a point $z\in Z$, i.e. $\d_Z$ is a Borel probability measure such that
for any Borel set $A\ss X$ one has $\d_z(A)=1\Liff z\in A$. Obviously, for any function $f\in C(X)$ one has
$$\int_Xfd\d_Z=\frac 1{|Z|}\sum\limits_{z\in Z}f(z).$$

\begin{Def} \label{st-hypap} In conditions of Definition \ref{hypap} let $\{Y_n\ |\ n\in\N\}$ be a sequence of
finite subsets of $X$. We say that the sequence $Y_n$ approximates the measure space $(X,\nu)$ if the sequence of measures
$\d_{Y_n}$ converges to the measure $\nu$ in the *-weak topology on the space $\M(X)$ of all Borel measures on $X$.
\end{Def}

The following proposition follows easily from well known theorems of functional analysis.

\begin{Prop} \label{ex-st-hypap}  In conditions of Definition \ref{st-hypap} suppose that every open ball in $X$ has the positive measure $\nu$ and
every set of the positive measure $\nu$ is infinite. Then for every set $A\ss X$ with $\nu(A)=1$ there exists a sequence $Y_n$ of finite subsets
of $X$ approximating the measure space $(X,\nu)$ such that $\all\,n\in\N\ Y_n\ss A$. $\Box$
\end{Prop}

Definition \ref{st-hypap} can be considered as a standard sequence version of Definition \ref{hypap} (1). This statement is based on the following

\begin{Prop} \label{equivalence} A sequence $Y_n\ss X$ approximates a measure space $(X,\nu)$ in the sense of Definition \ref{st-hypap} if and only if
for any $N\in\iN$ the set $Y_N$ is a hyperfinite approximation of the measure space $(X,\nu)$. \end{Prop}

\textbf{Proof.} $\Limpl$ Let $Y_n$ approximates $(X,\nu)$ and $N\in\iN$. Then for any $f\in C(X)$ one has
\begin{equation} \label{equivalence1}
\int f(\sp(x))d\mu_L\o\left(\frac 1{|Y_N|}\sum\limits_{y\in Y_N}\*f(y)\right)=\int\limits_Xfd\nu
\end{equation}
The first equality is due to $\*f$ is a lifting of $f\circ\sp$.
The second follows from Definition \ref{st-hypap} and from the nonstandard analysis definition of the limit of a sequence.
Now $\sp:Y_N\to X$ defines a measure $\nu'$ on $X$ that is the image of the Loeb measure of $Y$. Due to (\ref{equivalence1}) and the Riesz
representation theorem $\nu'=\nu$.

$\Longleftarrow$ Assume that $\sp\upharpoonright Y_N:Y_N\to\R$ is a measure preserving transformation for every $N\in\iN$. It is easy to see that for every function
$f\in C(X)$ the internal function $\*f\upharpoonright Y_N$ is a lifting of $f$. So,
$$
\o\left(\frac 1{|Y_N|}\sum\limits_{y\in Y_N}\*f(y)\right)=\int\limits_{Y_N}\o(\*f)d\mu_L=\int\limits_{Y_N}f\circ st\,d\mu_L=\int\limits_Xfd\nu.
$$
Thus, the equality (\ref{equivalence1}) holds for every $N\in\iN$ and by the nonstandard analysis definition of a limit one has
$\lim\limits_{n\to\infty}\int\limits_Xfd\d_{Y_n}=\int\limits_Xfd\nu$ $\Box$

The last two propositions imply Theorem \ref{cycle-appr}(1).

We use the same approach as above to formulate a sequence version of the notion of a hyperfinite approximation of a dynamical system.

\begin{Def} \label{seqap} Let $(X,\rho)$ be a compact metric space, $\nu$ be a Borel measure on
$X$, $\tau:X\to X$ be a measure preserving transformation of $X$, $\{Y_n\ss X\ |\ n\in\N\}$
be a sequence of finite approximation of the measure space $(X,\nu)$ in the sense of
Definition \ref{st-hypap} and $T_n:Y_n\to Y_n$ be a sequence of permutations of $Y_n$.
We say that a sequence $T_n$ is an approximating sequence of the transformation $\tau$ if for every $N\in\iN$ the internal permutation $T_N:Y_N\to Y_N$ is a h.a. of
$\tau$ in the sense of Definition \ref{hypap}. In this case we say that the sequence of finite dynamical systems $(Y_n,\mu_n, T_n)$ approximates the dynamical system
$(X,\nu,\tau)$. Here $\mu_n$ is a uniform probability measure on $Y_n$. 
\end{Def}

The reformulation of this definition in full generality in standard mathematical terms is practically unreadable. However, it is easy to reformulate it
for the case of an almost everywhere continuous transformation $\tau$. This case covers a lot of important applications.

Denote the set of all points of continuity of the map $\tau:X\to X$ by $D_{\tau}$.

\begin{Lm} \label{seqap1} Suppose that $\nu(D_{\tau})=1$ and let $Y\ss X$ be a h.a. of the measure space $(X,\nu)$. Then a permutation $T:Y\to Y$ is a h.a. of the
transformation $\tau$ if and only if for every positive $\e\in\R$ one has
\begin{equation} \label{seqap2}
\frac 1{|Y|}\left({|\{y\in Y\ |\ \rho(T(y),\*\tau(y))>\e\}|}\right)\approx 0.
\end{equation}\end{Lm}

\textbf{Proof} ($\Limpl$) Let $A=\{y\in Y\ |\ \o y\in D_{\tau}\}=\sp^{-1}(D_{\tau})$, $B=\{y\in Y\ |\ T(y)\approx\tau(\o y)\}$. Then, $\mu_L(A)=1$ since $Y$ is a h.a. of
the measure space $(X,\nu)$ and $\nu(D_{\tau})=1$. Since $T$ is a h.a. of $\tau$, one has $\mu_L(B)=1$. Thus, $\mu_L(A\cap B)=1$. Since $\tau$ is continuous on
$\D_{\tau}$ one has
\begin{equation}\label{seqap3}
\all\,x\in X\ \o x\in D_{\tau}\Limpl \*\tau(x)\approx\tau(\o x).
\end{equation}
So, $\all\, y\in A\cap B\ \*\tau(y)\approx\tau(\o y)$ and thus,  $\all\, y\in A\cap B\ \*\tau(y)\approx T(y)$. So, for every positive $\e\in\R$ one has
$\{y\in Y\ |\ \rho(T(y),\*\tau(y))>\e\}\ss Y\setminus(A\cap B)$. This proves (\ref{seqap2}).

$(\Longleftarrow)$ Suppose that (\ref{seqap2}) holds for every positive $\e\in\R$. Then obviously $\mu_L(\{y\in Y\ |T(y)\approx\*\tau(y)\})=1$. On the other hand,
by (\ref{seqap3}) one has $\mu_L(\{y\in Y\ |\ \*\tau(y)\approx\tau(\o y)\})=1$. Thus, $\mu_L(\{y\in Y\ |\ T(y)\approx\tau(\o y)\})=1$, i.e. $T$ is a h.a. of $\tau$
$\Box$

Lemma \ref{seqap1} implies immediately the following

\begin{Prop}[Standard version of Definition \ref{seqap}] \label{seqap4} In conditions of Definition \ref{seqap} and Lemma \ref{seqap1} the sequence of permutations
$T_n:Y_n\to Y_n$ is an approximating sequence of the transformation $\tau$ 
if and only if for every positive $\e\in\R$ one has
\begin{equation}\label{seqap5}
\lim\limits_{n\to\infty}\frac 1{|Y_n|}\left({|\{y\in Y_n\ |\ \rho(T_n(y),\tau(y))>\e\}|}\right)=0.\qquad \Box
\end{equation}
\end{Prop}

Now we are able to prove the sequence version of Theorem \ref{cycle-appr} (2).

\begin{Th}\label{exseqap} Let $(X,\rho)$ be a compact metric space and $\nu$ be a Borel measure on $X$ such that the measure space $(X,\nu)$ satisfies the conditions
of Theorem~\ref{cycle-appr}. Then for every measure preserving transformation $\tau:X\to X$ with $\nu(D_{\tau})=1$ there exist a sequence of finite sets
$Y_n\ss X$ and a sequence of permutations $T_n:Y_n\to Y_n$ such that the sequence of finite dynamical systems $(Y_n,T_n)$ approximates the dynamical system $(X,\nu,\tau)$
in the sense of Definition \ref{seqap}. Moreover, one can choose transitive permutations $T_n$.
\end{Th}

\textbf{Proof.} Let $Y_n\ss X$ be a sequence that approximates the measure space $(X,\nu)$ in the sense of Definition \ref{st-hypap}. Such sequence exists by
Proposition \ref{ex-st-hypap}. Then by Proposition \ref{equivalence} for any $N\in\iN$ the set $Y_N$ is a h.a. of the measure space $(X,\nu)$ in the sense of Definition
\ref{hypap}. By Theorem \ref{cycle-appr} there exists a (transitive) permutation $T_N:Y_N\to Y_N$ that is a h.a. of the transformation $\tau$.
By Lemma \ref{seqap1}, since $\nu(D_{\tau})=1$, this means that $(Y_N,T_N)$ satisfies (\ref{seqap2}) for every standard positive $\e$. In this proof the
letter $T$ maybe with lower indexes always denotes a (transitive) permutation.

For every numbers $n,m\in\N$ define the set
$$A_{n,m}=\left\{k\in\N\ \left| \right.\ \ex\, T:Y_k\to Y_k\
\left(\frac 1{|Y_k|}\cdot\left| \{y\in Y_k\ |\ \rho(T(y),\tau(y))>\frac 1n\}\right|<\frac 1m\right)\right\}.$$

Since $\all\,N\in\iN\ N\in\*A_{n,m}$, there exists a standard function $N(n,m)$ such that $\all\,k>N(n,m)\ k\in A_{n,m}$. By the definition of sets $A_{m,n}$, there exists
a standard function $T(k,n,m):Y_k\to Y_k$ with the domain $\{(n,m,k)\in\N^3\ |\ k>N(n,m)\}$ such that
$$
\frac 1{|Y_k|}\cdot\left| \left\{y\in Y_k\ |\ \rho(T_k(y),\tau(y))>\frac 1n\right\}\right|<\frac 1m.
$$

Now it is easy to see that if $r=N(n,n)+n$, then the sequence $(Y_r,T_r)$ satisfies the conditions of Proposition \ref{seqap4} $\Box$

This proof is based on NSA. The purely standard proof of this theorem seems to be much more difficult.

\bigskip

\section{Proofs of Theorems \ref{NSBET}, \ref{cycle-appr}(2) and \ref{un-erg}.}

\bigskip

i) \emph{(Proof of Theorem \ref{NSBET})}). We notice first, that the following proposition follows immediately from
Theorem \ref{BET} applied to the external dynamical system $(Y,\mu_L,T)$ and Theorem \ref{sequences}.

\begin{Prop}\label{WeakNSBET} In conditions of Theorem \ref{ErgMeanStab} for any $y\in Y$
there exists an initial segment $I\ss\iN$ such that
$$\all\,L,K\in I\ (A_K(F,T,y)\approx A_L(F,T,y)\approx\lim\limits_{n\to\infty}A_n(\o F,T,y)).$$
\end{Prop}

Proposition \ref{WeakNSBET} is a weaker version of Theorem \ref{NSBET}, since it differs of this theorem only in the order of quantifiers "for all $y\in T$
and  "there exists an initial segment $I\ss\iN$".

To prove Theorem \ref{NSBET} first it is necessary to prove

\begin{Th} \label{a.e.-convergence} Let $f_n:Y\to\R,\ n\in\N$ be a sequence
of $\mu_L$ measurable functions on $Y$, and $F_n:Y\to\*\R,\ n\in\*N$ be an
\textbf{internal} sequence such that $\all\,n\in\N,\ F_n$ is a
lifting of $f_n$. Then $f_n$ converges to a measurable function
$f$ $\mu_L$-almost everywhere if and only if there exists $K\in\iN$ such that
 $\mu_L$-almost everywhere $\all\, N\in\iN,\ N<K\Longrightarrow F_N(x)\approx
F(x)$, where $F$ is a lifting of $f$.
\end{Th}

\textbf{Proof} of Theorem \ref{a.e.-convergence}.

($\Longrightarrow$) Let $f_n$converges to $f$ a.e.  By Egoroff's
theorem
$$
\all k\in\N\,\exists B_k\subseteq Y\, (\mu_L(B_k)\geq 1-\frac
1k\land f_n(x)=\o F_n(x)\ \mbox{converges uniformly on}\
B_k).
$$
WLOG we may assume that $B_k$ is internal, $\frac{|B_k|}{|Y|}\geq
1-\frac 1k$, and $\all\, n,k\in\N\,\all x\in B_k\ F_n(x)\approx f_n(x)$
and $F(x)\approx f(x)$. Then
$$
\exists\st\f_k:\N\to\N\,\all^{st}
r\,\all\st m>\f_k(r)\max\limits_{x\in
B_k}\,|F_m(x)-F(x)|<\frac 1r.
$$

Consider the internal set
$$
C^k_r=\{N\in\*\N\ |\ \all\,m \,(N>m>\f_k(r)\Limpl\all x\in B_k\
|F_m(x)-F(x)|<\frac 1r)\}
$$
The previous statement shows that $C^k_r$ contains all standard
numbers that are greater that $\f_k(r)$. Thus, there exists infinite $N^k_r\in
C^k_r$. By Theorem \ref{spill} (5)  $\exists K\in\iN\, \all\st k,r\,
K<N^k_r$. Obviously, this $K$ satisfies Theorem \ref{a.e.-convergence}

($\Longleftarrow$) Let $B=\{x\in Y\ |\ \all N\in\iN\;\;\; N\leq K\Limpl
F_N(x)\approx F(x)\},\ A_n=\{x\in Y\ |\ f_n(x)\approx F_n(x)\},\
n\in\N,\ A=\{x\in Y\ |\ F(x)\approx f(x)\}, \ C=B\cap
A\cap\bigcap\limits_{n\in\N}A_n.$

By conditions of the theorem $\mu_L(C)=1$. Fix an arbitrary $x\in C$, and an arbitrary
$r\in\N$. The internal set $D=\{n\in\*\N\ \ |\ |F_n(x)-F(x)|\leq \frac 1r\}$ contains
all infinite numbers that are less or equal to $K$. So $\exists n_0\in\N\, \all n>n_0\
|F_n(x)-F(x)|\leq \frac 1r\}$. Since $F_n(x)\approx f_n(x)$, the same holds for
$f_n(x)$ and $\o F(x)$. Thus, $f_n$ converges to $f=\o F$ a.e.     $\Box$.

Now we are able to complete the proof of Theorem \ref{NSBET}. In conditions of Theorem
\ref{NSBET} let $f=\o F$ and $f_n(x)=A_n(f,T,x),\ n\in\N$ and
the internal sequence $F_n(x)=A_n(F,T,x),\ n\in\*\N$. Then $f\in
L_1(\mu_L)$ and $F_n$ is an $S$-integrable lifting of $f_n$ for
all $n\in\N$. By Theorem \ref{BET} $f_n$ converges to an
integrable function $\hat f$ a.e. Let $\hat F$ be an
$S$-integrable lifting of $\hat f$. Then by Theorem
\ref{a.e.-convergence} there exists $K\in\iN$ such that
$\mu_L$-almost surely $\all N\in\iN\ N<K\Longrightarrow
F_N(X)\approx \hat F(x)$ $\Box$.

ii) (\emph{Proof of Theorem \ref{cycle-appr} (2)}.) It is well-known (see e.g. \cite{KSF,Br}) the measure space $(X,\nu)$ is a Lebesgue space, i.e. it
is isomorphic modulo measure $0$ to the measure space $([0,1],dx)$, where $dx$ is the standard Lebesgue measure.
This means that there exist a set $B\ss X$ a set $C\ss [0,1]$ and a bijective map $\psi:B\to C$
such that $dx(C)=\nu(B)=1$ and the maps $\psi,\psi^{-1}$ are measure preserving.

\begin{Lm} \label{hyp-iso} In conditions of the previous paragraph let $Y$ be a h.a. of $(X,\nu))$. Then for every set $D\ss [0,1]$
with $dx(D)=1$ there exists a bijective lifting $G:Y\to\*[0,1]$ of the map $\psi$ such that
\begin{enumerate}
\item $Z=G(Y)\ss \*D$;
\item $Z$ is a h.a. of $([0,1],dx)$.
\item $G^{-1}:Z\to \*X$ is a lifting of $\psi^{-1}$.
\end{enumerate}
\end{Lm}

\textbf{Proof}. Let $F:Y\to\*[0,1]$ be a lifting of $\psi$. Let $\s=\frac12\min\{\rho(u,v)\ |\ u,v\in F(Y),\,u\neq v\}$.
Then $0<\s\approx 0$ and $\all\, u\in F(Y)\ B_{\s}(u)\cap F(Y)=\{u\}$. Since $\*\nu (B_{\s}(u))>0$ and $dx(D)=1$ the set $B_{\s}(u)\cap\*D$ contains infinitely
many points, and thus, there exists an internal set $E_u\ss B_{\s}(u)\cap\*D$ such that $|E_u|=|F^{-1}(u)|$. Establishing bijection between $F^{-1}(u)$ and $E_u$ for every
$u\in F(Y)$, we obtain the bijection $G:Y\to Z\ss\*D$ that is a lifting of $\psi$. Notice that since $G$ and $G^{-1}$ are bijections they are measure preserving
maps between measure spaces $(Y,\mu_L^Y)$ and $(Z,\mu_L^Z)$.

To prove the second property of the set $Z$, one needs to show that $st_{[01]}\upharpoonright Z:Z\to [0,1]$ is a measure preserving map, i.e. that for every measurable set
$A\ss[0,1]$ one has
\begin{equation} \label{hyp-iso1}
\mu_L^Z(st^{-1}_{[0,1]}(A)\cap Z)=dx(A)
\end{equation}
One has
\begin{equation}\label{hyp-iso2}
\mu_L^Z(st^{-1}_{[0,1]}(A)\cap Z)=\mu_L^Y(\{y\in Y\ |\ G(y)\in st^{-1}_{[0,1]}(A)\})=\mu_L^Y(\{y\in Y\ |\ \o G(y)\in A\}).
\end{equation}
Since $G$ is a lifting of $\psi$ on has $\o G(y)=\psi(st_X(y))$ for $\mu_L^Y$-almost all $y$. Thus,
\begin{equation}\label{hyp-iso3}
\mu_L^Y(\{y\in Y\ |\ \o G(y)\in A\})=\mu_L^Y(\{y\in Y\ |\ \psi(st_X(y))\in A\})=\nu(\psi^{-1}(A))=dx(A),
\end{equation}
since $st_X\upharpoonright Y:Y\to X$ and $\psi:B\to C\ss [0,1]$ are measure preserving maps.
The equality (\ref{hyp-iso1}) follows from the equalities (\ref{hyp-iso2}) and (\ref{hyp-iso3}).

To prove the third property of the set $Z$ it is enough to show that $st_Y(G^{-1}(z))=\psi^{-1}(st_{[0,1]}(z))$ for $\mu_L^Z$-almost all $z\in Z$. Since $\psi$ is a bijection,
the last equality is equivalent to the equality $\psi (st_Y(G^{-1}(z)))=st_{[0,1]}(z)$, which follows from the following sequence of equalities  that hold for $\mu_L^Z$-almost
all $z\in Z$:
$$\psi(st_Y(G^{-1}(z))=st_{[0,1]}(G(G^{-1}(z)))=st_{[0,1]}(z). \Box$$

The proof of Theorem~\ref{cycle-appr}(2) is divided in six parts I -- VI.

I. Here we prove the existence of a h.a. $(Y,\mu_L,T)$ of the
dynamical system $([0,1],dx,\tau)$.
Let $Y$ be an arbitrary h.a. of the measure space $([0,1],dx)$.  Let $F:Y\to\*[0,1]$ be a lifting of $\tau$.
First we prove  the following statement.

(A) \emph{For every standard $\d>0$ there exists a permutation
$T_{\d}:Y\to Y$ such that}
$$
\frac{|\{y\in Y\ |\ |F(y)-T_{\d}(y)|<\d\}|}M\approx 1.
$$
We deduce (A) from the Marriage Lemma.
Fix a standard $\d>0$ and for every $y\in Y$ set
$S(y)=\*(F(y)-\d,F(y)+\d)\cap Y$. Let $I$ be an arbitrary internal
subset of $Y$. Set $S(I)=\bigcup\limits_{y\in I}S(y)$ and
$B(I)=\bigcup\limits_{y\in I}\*(F(y)-\d,F(y)+\d)$. So,
$S(I)=B(I)\cap Y$. The internal set $B(I)$ can be represented as a
union of a hyperfinite family of disjoint intervals. Since the
length of each of these intervals is not less than $2\d$, their
number is actually finite. Let
$B(I)=\bigcup\limits_{i=1}^n(\xi_i,\eta_i)$, where intervals
$(\xi_i,\eta_i)$ are pairwise disjoint and $n$ is standard.

Consider the standard set
$C=\bigcup\limits_{i=1}^n(\o\xi_i,\o\eta_i)$. Then $dx(C)\mu_L(\sp^{-1}(\*C))$. Obviously, $\sp^{-1}(\*C)\Delta
B(I)\ss\bigcup\limits_{i=1}^n(\Mu(\o\xi_i)\cup \Mu(\o\eta_i))=\Mu(\partial C)$,
where the monad of a number $a\in[0,1]$ is denoted by $\Mu(a)$. Since the Loeb measure
of the monad of any number is equal to $0$ and so, $\Mu(\partial C)=0$, one has
$dx(C)=\o\left(\frac{|S(I)|}M\right)$. Substituting $[0,1]$ for $X$, $dx$ for $\nu$ and $\tau$ for $\psi$ in (\ref{hyp-iso3}) obtain
$dx(C)=dx(\tau^{-1}(C))=\mu_L(F^{-1}(\sp^{-1}(C)))$. Since $I\setminus
F^{-1}(\sp^{-1}(C))\ss\Mu(\partial C)$, one has $\o\left(\frac{|I|}M\right)\leq
\o\left(\frac{|S(I)|}M\right)$. This means that if
$r_I=\max\{0,|I|-|S(I)|\}$, then $\frac{r_I}M\approx 0$. Let
$r=\max\{r_I\ | \ I\ss Y \}$. Fix an arbitrary set
$Z\ss\*[0,1]\setminus Y$ such that $|Z|=r$. For every $y\in Y$
set $S'(y)=S(y)\cup Z$ and for an arbitrary $I\ss Y$ set
$S'(I)=\bigcup\limits_{y\in Y}S'(y)$. Then $S'(I)=S(I)\cup Z$,
$|S'(I)|=|S(I)|+r\geq |I|$, since $|I|-|S(I)|=r_I\leq r$. By the
Marriage Lemma there exists an injective map $\theta:Y\to
S'(Y)=Y\cup Z$ such that $\all\,y\ \theta(y)\in S'(y)$. Obviously
$|\theta^{-1}(Z)|=|Y\setminus\theta(Y)|\leq r$. So, there exists a
bijective map $\l:\theta^{-1}(Z)\to Y\setminus\theta(Y)$. Define
$T_{\d}:Y\to Y$ by the formula
$$
T_{\d}(y)=\left\{\begin{array}{ll} \theta(y),\ &y\in
Y\setminus\theta^{-1}(Z) \\ \l(y),\ &
y\in\theta^{-1}(Z)\end{array}\right.
$$
Notice that $\frac{|\theta^{-1}(Z)|)}M\leq \frac rM$. By
construction of $T_{\d}$ one has $\all\,y\in
Y\setminus\theta^{-1}(Z)\ |T_{\d}(y)-\tau(y)|<\d$.
Since $\mu_L(\theta^{-1}(Z))\leq\frac rM\approx 0$, the statement
(A) is proved.

Let $\cs(Y)$ be the set of all internal permutations of $Y$. Consider
the external function $f:\N\to \cs(Y)$ such that $f(n)=T_\frac
1n$. By the Saturation Principle the function $f$ can be extended to an internal function $\bar
f:\{0,\dots,K\}\to \cs(Y)$ for some $K\in\iN$. Internal function
$g(n)=\frac{|\{y\in Y\ |\ |F(y)-\bar f(y)|\geq\frac 1n\}|}M$ assumes
only infinitesimal values for all standard $n$. By Robinson's
Lemma there exists $L\in\iN$ such that $g(L)\approx 0$. set
$T=\bar f(L)$. Then $\mu_L(\{y\in Y\ |\ T(y)\approx F(y)\})=1$.
Since $F$ is a lifting of $\tau$, the same is true also for
$T(y)$. This proves I.

We have to prove now that a h.a. $T$ of $\tau$ can be chosen as
a cycle of maximal length.

II. Fix a permutation $T:Y\to Y$ that is a h.a. of $\tau$ and represent it
by a product of pairwise disjoint cycles, including the cycles of length 1 (fix points):
\begin{equation}\label{II1}
T=(y_{11}...y_{1n_1})(y_{21}...y_{2n_2})...(y_{b1}...y_{bn_b}),
\end{equation}
where $y_{ij}\in Y$ is the $j$-th element in the $i$-th cycle and $b$ is the number of cycles.
So,
\begin{equation}\label{II2}
\sum\limits_{i=1}^bn_i=M=|Y|.
\end{equation}
We assume also that $n_1\geq n_2\geq\dots\geq n_b.$
Consider the cycle
\begin{equation}\label{II3}
C=(y_{11}...y_{1n_1}y_{21}...y_{2n_2}...y_{b1}...y_{bn_b})
\end{equation}
By (\ref{II2}) $C$ is a cycle of length $M$, i.e. a transitive permutation.

Set $B=\{y\in Y\ |\ C(y)\neq T(y)\}$.
\begin{equation}\label{II4}
|B|=b=\sum\limits_{n=1}^Ma_n,
\end{equation}
where $a_n$ is the number of cycles of length $n$.

III.  Recall that a point $x\in [0,1]$ is said to be an $n$-periodic point of the transformation $\tau$
if its orbit under this transformation consists of $n$-points:
$x,\tau x,\dots,\tau^{n-1}x$. A point $x$ is said to be $\tau$-periodic if it is
$n$-periodic for some $n$.  The transformation $\tau$ is said to
be aperiodic if the set of periodic
points has measure zero. It is well-known that every measure
preserving automorphism $\tau$ of a Lebesgue space $X$ defines the
partition of this space by $\tau$-invariant Lebesgue subspaces of
aperiodic and $n$-periodic points. So, it is enough to prove our
statement for the case of aperiodic transformation $\tau$ and for
the case of $n$-periodic transformation $\tau$.

Suppose that the transformation $\tau$ is aperiodic. Let us
prove that under this assumption the cycle $C$ defined in
the part II is a h.a. of $\tau$.

Let $P_n(T)\ss Y$ be the set of all $n$-periodic points of $T$ and let $P_n(\tau)\ss X$ be the
set of all $n$-periodic points of $\tau$.
Since $T$ is a lifting of $\tau$
it is easy to that for every standard $k$ the following relations
\begin{equation}\label{IV1}
T(y)\approx\tau(\o y),\dots,T^{k}(y)\approx\tau^{k}{(\o y)}
\end{equation}
hold $\mu_L$-a.e. on $Y$.
So, for every standard $n$ $P_n(T)\ss st^{-1}(P_n(\tau))$ up to a set of the Loeb measure zero. Since
$dx(P_n(\tau))=0$, one has $\frac 1M|P_n(T)|\approx 0$. Obviously, $|P_n(T)|=na_n$.
Thus, for every standard $n$ one has $\frac 1M\cdot
a_n\approx 0$.
By the Robinson's Lemma there exists an infinite $N$ such that $
\frac 1M\sum\limits_{n=1}^N a_n\approx 0.$. Obviously $
M\geq\sum\limits_{n=N+1}^Ma_n\cdot n\geq
(N+1)\sum\limits_{n=N+1}^M a_n.$ So, $\frac
1M\sum\limits_{n=N+1}^M a_n\leq\frac 1{N+1}\approx 0$ and $\frac
1M\cdot |B|=\frac 1M\sum\limits_{n=1}^Ma_n\approx 0. $ Thus,
$\mu_L(B)=0$, $C(y)=T(y)$ $\mu_L$-a.e. and $C$ approximates $\tau$.

IV. Suppose now that $\tau$ is $n$-periodic. We prove first that a
h.a. $T$ of $\tau$ also can be chosen to be $n$-periodic.
The relations (\ref{IV1}) imply that for almost every point $y\in Y$ if
$y$ has a standard period with respect to $T$, then this period is a multiple of $n$. Indeed, if
$y$ satisfies (\ref{IV1}), and its standard period is $nq+r$ for
$0<r<n$, then $\o y=\o T^{nq+r}(y)=\tau^{nq+r}(\o y)=\tau^{r}(\o
y)$, which is impossible since $\tau$ is $n$-periodic. By
Saturation Principle, there exist an internal set $I\ss Y$ such
that $\mu_L(I)=1$ and a number $N\in\iN$ such that for every point
$y\in I$, whose period is less, than $N$, this period is a
multiple of $n$.

Consider the representation (\ref{II1}) of $T$ and set $n_i=nq_i+r_i,\ r_i<n$ for each $i\leq b$.
Let $Y'\ss Y$ be the set obtained by deleting from $Y$ the last $r_i$ elements of the i-th cycle
for each $i\leq b$.
The set $Y'$ has the Loeb measure equal to
$1$. Indeed, all the deleted elements either belong to the set
$Y\setminus I$, whose measure is $0$, or to a cycle whose length
is greater, than $N$. The number of these cycles does not exceed
$\frac MN$ and the number of deleted points in each such cycle is
less, than $n$. So the Loeb measure of the set of these points is
also equal to $0$.
Since $\mu_L(Y')=1$ the pair $(Y',\sp)$ is a h.a. of $[0,1]$.
The construction of $Y'$ defines also the permutation $T':Y'\to Y'$
such that
\begin{equation}\label{V1}
T'=(y_{11}...y_{1\,n\cdot q_1})(y_{21}...y_{2\,n\cdot q_2})...(y_{b1}...y_{b\,n\cdot q_b}).
\end{equation}
Notice, that actually the number of cycles in $T'$ may be less, than $b$, since in case of $q_i=0$ the $i$-th cycle is empty.
However, the dynamical system $(Y',\mu_L,T')$ is a h.a. of the dynamical system $(X,\nu,\tau)$.
Indeed, let $D=\{y\in Y'\ |\ T(y)\neq T'(y)\}$. Then $D\ss\{y\in Y'\ |\ T(y)\in Y\setminus Y'\}\ss T^{-1}(Y\setminus Y')$.
Thus, $\mu_L(D)\leq \mu_L(Y\setminus Y')=0$
To obtain an $n$-periodic h.a. of $\tau$ it is enough to split each cycle in the representation (\ref{V1}) in cycles of length $n$.
Indeed, let the obtained cycle be
$$
T''=(z_1,...,z_n)(z_{n+1},...,z_{2\cdot n})\dots(z_{(i-1)\cdot n+1},...,z_{i\cdot n})\dots (z_{(K-1)\cdot n+1},...,z_{K\cdot n}),
$$
where $K=|Y'|/n$.
It is easy to see that $T''(y)\neq T'(y)$, only for the points
$z_{i\cdot n}$. Notice, that $\mu_L(\{z_{i\cdot n}\ |\ i\leq K\})=\frac 1n>0$
However, due to (\ref{IV1}) and the $n$-periodicity of $\tau$, for almost all of these points one has
$$T'(z_{i\cdot n})=z_{in+1}\approx\tau(\o z_{i\cdot n+1})=\tau^n(\o z_{(i-1)\cdot n+1})=\o z_{(i-1)\cdot n+1}.$$
At the same time $T''(z_{i\cdot n})=z_{(i-1)\cdot n+1}$ by the
definition. Thus, $T''(y)\approx T'(y)$ for almost all $y$.

V. To complete the proof of the theorem for $X=[0,1]$ we need to consider the case when all
orbits of $T$ have the same standard period $n$. In this case
$M=N\cdot n$.

It is easy to see that there exists a selector $I\subset Y$ 
(subset that intersect each orbit of $T$ by a single point) that is dense in
$\*[0,1]$, i.e. the monad $M(I)=\*[0,1]$. It is enough to show the
existence of a selector that intersects every interval with
rational endpoints. Obviously, for every finite set $A$ of such
intervals, there exists a selector that intersects each
interval from $A$. The existence of a dense selector follows from
the Saturation Principle.

Let $I=\{y_1<y_2<\dots<y_N\}$ be a dense selector. Here $<$ is the
order in $\*[0,1]$. Due to the density of $I$ in $\*[0,1]$ for
every $k<N$ one has $y_k\approx y_{k+1}$. Obviously, the transformation $T$
can be represented by a product of pairwise disjoint cycles as follows:
$$
T=(y_1,...,T^{n-1}y_1)(y_2,...,T^{n-1}y_2)\dots(y_N,...,T^{n-1}y_N)
$$
Consider the following cycle $S$ of the length $M$:
$$
S=(y_1,...,T^{n-1}y_1y_2,...,T^{n-1}y_2\dots y_N,...,T^{n-1}y_N)
$$
Since for every $k\leq N$ holds $T^n(y_k)=y_k$, one has
$$\o S(T^{n-1}(y_k))=\o y_{k+1}=\o y_k =\o T^n(y_k)=\o
T(T^{n-1}(y_k))=\tau(\o T^{n-1}(y_k))$$ for almost all $k$. Thus,
$\o S(y)=\tau(\o y)$ for almost all $y$ and the cycle $S$ is a h.a. of $\tau$.

We proved actually that for every h.a. $Y$ of $([0,1]),dx)$ there exists an internal set $Y'\ss Y$ with $\mu_L(Y')=1$
and a permutation $T':Y'\to Y'$ such that the hyperfinite dynamical system $(Y',\mu_L,T')$ is a h.a. of the dynamical system
$(X,\nu,\tau)$ and $T'$ is a transitive permutation of $Y'$ (see Part IV of this proof). To obtain a transitive h.a. $T:Y\to Y$ of $\tau$, set
$T'=(z_1,...,z_{|Y'|})$ and $Y\setminus Y'=\{u_1,...,u_{|Y\setminus Y'|}\}$ and consider the cycle of the length $|Y|$
$$
T:(z_1,...,z_{|Y'|},u_1,...,u_{|Y\setminus Y'|})
$$
Since $\mu_L(\{y\in Y\ |\ T'(y)\neq T(y)\})=0$ the transformation $T$ is h.a. of $\tau$.

VI. The statement of the theorem for the case of an arbitrary dynamical system 
$(X,\nu,\tau)$, satisfying the conditions, follows immediately from Lemma \ref{hyp-iso}.
Indeed, let a set $B\ss X$, a set $C\ss [0,1]$, a bijective map $\psi:B\to C$ and a bijective lifting $G:Y\to\*[0,1]$ of $\psi$
satisfy the conditions of Lemma \ref{hyp-iso}. Then $\lambda=\psi\tau\psi^{-1}:[0,1]\to [0,1]$ is a measure preserving
transformation. Fix an arbitrary h.a. $Y$ of the measure space $(X,\nu)$. Then by Lemma \ref{hyp-iso} the hyperfinite
set $Z=G(Y)$ is a h.a. of $([0,1],dx)$. By the results proved in the parts I-V, there exists a permutation $S:Z\to Z$
that is a h.a. of $\lambda$. Then it is easy to see that the permutation $T=G^{-1}SG:Y\to Y$ is a h.a. of $\tau$.
Obviously, if $S$ is a transitive permutation, then $T$ is a transitive permutation as well. $\Box$

\bigskip

iii) (\emph{Proof of Theorem \ref{un-erg}.}) Let $(X,\rho)$ be a compact metric space. Consider a hyperfinite set $Y\ss \*X$. This set defines a
Borel measure $\nu_Y$ on $X$ by the formula $\nu_Y(K)=\mu_L(st^{-1}(K)\cap Y)$. Obviously $Y$ is a h.a.
of the measure space $(X,\nu_Y)$. Let $T:Y\to Y$ be an internal permutation that is
$S$-continuous on $A$
for some (not necessary internal) set
$A\ss Y$ with $\mu_L(A)=1$, i.e.
\begin{equation} \label{measurable}
\all\, a_1,a_2\in A\ (a_1\approx a_2\Limpl T(a_1)\approx T(a_2)).
\end{equation}
Notice that since $\sp^{-1}(st(A))\supseteq A$ and $\mu_L(A)=1$,
the set $\sp(A)\ss X$ is a measurable set w.r.t. the completion of the measure $\nu_Y$, which we denote by $\nu_Y$ also,
and $\nu_Y(\sp(A))=1$.

Define a map $\tau_T:X\to X$ such that $\tau_T(\sp(y))=\sp(T(y))$ for $y\in A$ and $\tau_T\upharpoonright{X\setminus\sp(A)}$
is an arbitrary measurable permutation of the set $X\setminus\sp(A)$.
\begin{Prop} \label{tauT}
The map $\tau_T$ preserves the measure $\nu_Y$.
\end{Prop}
\textbf{Proof}. Replacing, if necessary, $A$ by $\bigcap\limits_{n\in\N}T^n(A)$ we may assume that $A$ is invariant for permutation $T$.
Then, obviously, $\sp(A)$ is invariant for $\tau_T$.

Consider a closed set $B\ss X$. We have to prove that $\nu_Y(\tau^{-1}_T(B))=\nu_Y(B)$. One has
$$\nu_Y(\tau_T^{-1}(B))=\nu_Y(\tau^{-1}_T(B)\cap\sp(A)).$$
It is easy to check that
$$\tau_T^{-1}(B)\cap\sp(A)=\sp(T^{-1}(\sp^{-1}(B))\cap A).$$
Thus,
$$\nu_Y(\tau_T^{-1}(B))=\mu_L\left(\sp^{-1}(\sp(T^{-1}(\sp^{-1}(B))\cap A))\right)=\mu_L\left(\left(\sp^{-1}(\sp(T^{-1}(\sp^{-1}(B))\cap A))\right)\cap A\right)$$
Using (\ref{measurable}) and the $T$-invariance of it is easy to check, that
$$\sp^{-1}(\sp(T^{-1}(\sp^{-1}(B))\cap A))\cap A=T^{-1}(\sp^{-1}(B))\cap A).$$
So,
$$\nu_Y(\tau_T^{-1}(B))=\mu_L(T^{-1}(\sp^{-1}(B))\cap A))=\mu_L(T^{-1}(\sp^{-1}(B))=\mu_L(\sp^{-1}(B))=\nu_Y(B). $$
In the last chain of equalities we used the facts that $\mu_L(A)=1$ and that $T$ being a permutation preserves the Loeb measure. $\Box$

\begin{Prop}\label{functional}
1) In conditions of Proposition \ref{tauT} for any $a>0$ and for any $y\in Y$ the following positive functional $l_a(\cdot,T,y)$ on $C(X)$ is defined:
$l_a(f,T,y)=\o A_K(\*f,T,y)$, where $\o\left(\frac KM\right)=a$ and $f\in C(X)$

2) If $\all\,K,L\in\iN\ (\frac KM\approx\frac LM\approx 0\Limpl A_K(\*f,T,y)\approx A_L(\*f,T,y)$), then $l_0(f,T,y)$ is defined by the same formula as in 1).
In this case $l_0(f,T,y)=\wt f(y)$

3) If $T:Y\to Y$ is $S$-continuous, then the functional $l_0(\cdot,T,y)$ is $\tau_T$-invariant for all $y\in Y$.
\end{Prop}

\textbf{Proof}. The correctness of the definition in 1) follows from Theorem \ref{ErgMeanStab}. The statement 2) follows from Proposition \ref{WeakNSBET}. To prove
statement 3) notice that if $T$ is $S$-continuous on $Y$, then $\tau_T$ is continuous on $X$ and, thus, $\*(f\circ\tau_T)\upharpoonright Y$ is a lifting of $f\circ\tau_T$. So,
$$\all\, y\ \all\, K\in\*\N\ \*(f\circ\tau_T) (T^K(y))\approx f(\tau_T(\o T^K(y)))=f(\o T^{K+1}(y))\approx\* f(T^{K+1}y)$$
These equivalences allows to prove that $A_K(\*(f\circ\tau),T,y)\approx A_K(\* f,T,y)$. $\Box$

Now we can complete the proof of Theorem \ref{un-erg}

Let $y\in Y$ satisfy conditions of the theorem. For a number $K\in\*\N$
we denote the initial $K$-segment of the $T$-orbit of $Y$ by $S(K,y)$. Then
for any $K\in\iN$ one has $\sp(S(K,y))=Y$, since the closed set $\sp(S(K,y))$ contains the $\tau$-orbit of $\sp(y)$. Let $K$ be a $T$-period of $y$. Then $K\in\iN$. Otherwise, the $\tau$-orbit of $\sp(y)$ would be finite, while we assume $X$ to be infinite. It is easy to see that it is enough to prove the theorem
for every $N\in\iN$ such that $N\leq K$. Under this assumption all elements of the set $Y_1=\{y, Ty, ... , T^{N-1}y\}$ are distinct. Since $\sp(Y_1)=X$, the set $Y_1$
defines the Borel measure $\nu_{Y_1}$ on $X$. Let $T_1:Y_1\to Y_1$ be the permutation of $Y_1$ that differs from $T$ only for one element $T^{N-1}y$: $T_1(T^{N-1}y)=y$.
Set $A=Y_1\setminus\{T^{N-1}y\}$. Then $X$, $\tau$, $Y_1$, $T_1$, and $A$ satisfy conditions of Proposition \ref{tauT}: $\mu_L(A)=1$,
$\all\, z\in A\ \sp(T_1z) =\tau(\sp(z))$, i.e. $\tau_{T_1}=\tau$ and $T_1$ is $S$-continuous on $A$, since $\tau$ is a continuous map. By Proposition \ref{tauT}
the measure $\nu_{Y_1}$ is $\tau$-invariant. Thus, $\nu_{Y_1}=\nu$ due to the unique ergodicity of the map $\tau$. If $f\in C(X)$, then obviously $\*f\upharpoonright Y_1$ is an $S$-integrable lifting of $f$. This proves the equality (\ref{un-erg-1}). $\Box$

Gordon at Eastern Illinois University 600 Lincoln Avenue
Charleston, IL 61920-3099 USA, email: cfyig@eiu.edu.

Glebsky at IICO-UASLP
AvKarakorum 1470
Lomas 4ta Session
SanLuis Potosi SLP 7820 Mexico, email: glebsky@cactus.iico.uaslp.mx

Henson at University of Illinois at Urbana-Champaign 1409 W. Green Street
Urbana, Illinois 61801-2975 USA, email: henson@math.iuc.edu

\end{document}